\documentclass[10pt]{article}

\usepackage{pifont}
\usepackage{amsmath}
\usepackage{cases}
\usepackage{amsfonts}
\usepackage{latexsym}
\usepackage{amssymb}
\usepackage{stmaryrd}
\usepackage{overpic}
\usepackage{color}
\usepackage{abstract}
\usepackage{indentfirst}
\usepackage{booktabs}
\usepackage{url}
\usepackage{verbatim}
\usepackage{multirow} 
\usepackage{graphicx}  
\usepackage{float}  
\usepackage{subfigure} 
\usepackage{algorithm}
\usepackage{algorithmic}
\newcommand{\R}{\mathbb{R}}
\newcommand{\A}{\mathrm{A}}

\newcommand{\C}{\mathcal{C}}
\newcommand{\bu}{\boldsymbol{u}}
\newcommand{\bx}{\boldsymbol{x}}
\newtheorem{theorem}{Theorem}[section]

\newtheorem{example}[theorem]{Example}

\newtheorem{remark}[theorem]{Remark}

\begin{document}
	
	\begin{center}
		\large\bf  Overlapping Schwarz Preconditioners for Randomized \\ Neural Networks with Domain Decomposition 
	\end{center}
	
	\vspace{2mm}
	
	\begin{center}
	{\large Yong Shang}\footnote{School of Mathematics and Statistics, Xi'an Jiaotong University, Xi'an, Shaanxi 710049, P.R. China. E-mail: {\tt fsy2503@stu.xjtu.edu.cn}},\quad
	{\large Alexander Heinlein}\footnote{Delft Institute of Applied Mathematics, Faculty of Electrical Engineering, Mathematics and Computer Science, Delft University of Technology, Mekelweg 4, 2628 CD Delft, Netherlands. E-mail: {\tt a.heinlein@tudelft.nl}},\quad
	{\large Siddhartha Mishra}\footnote{Seminar for Applied Mathematics, D-MATH, and ETH AI Center, ETH Z\"urich, Rämistrasse 101, 8092 Zürich, Switzerland. E-mail: {\tt smishra@sam.math.ethz.ch}},\quad {\rm and}	\quad
	{\large Fei Wang}\footnote{School of Mathematics and Statistics, Xi'an Jiaotong University, Xi'an, Shaanxi 710049, P.R. China. The work of this author was partially supported by the National Natural Science Foundation of China (Grant No.\ 92470115). Email: {\tt feiwang.xjtu@xjtu.edu.cn}}

	\end{center}

	\vspace{2mm}
	
	\begin{quote} 
		\noindent{}{\bf Abstract}: Randomized neural networks (RaNNs), in which hidden layers remain fixed after random initialization, provide an efficient alternative for parameter optimization compared to fully parameterized networks. In this paper, RaNNs are integrated with overlapping Schwarz domain decomposition in two (main) ways: first, to formulate the least-squares problem with localized basis functions, and second, to construct overlapping preconditioners for the resulting linear systems. In particular, neural networks are initialized randomly in each subdomain based on a uniform distribution and linked through a partition of unity, forming a global solution that approximates the solution of the partial differential equation. Boundary conditions are enforced through a constraining operator,  eliminating the need for a penalty term to handle them.  Principal component analysis (PCA) is employed to reduce the number of basis functions in each subdomain, yielding a linear system with a lower condition number. By constructing additive and restricted additive Schwarz preconditioners, the least-squares problem is solved efficiently using the Conjugate Gradient (CG) and Generalized Minimal Residual (GMRES) methods, respectively. Our numerical results demonstrate that the proposed approach significantly reduces computational time for multi-scale and time-dependent problems. Additionally, a three-dimensional problem is presented to demonstrate the efficiency of using the CG method with an AS preconditioner, compared to an QR decomposition, in solving the least-squares problem.
		
		
	\end{quote}
	
	\vspace{2mm}

	\section{Introduction}

	Domain decomposition methods (DDMs) are based on dividing the computational domain into overlapping or non-overlapping subdomains, breaking the problem into smaller subproblems. This allows local subproblems to be solved efficiently on multiple processors. DDMs greatly improve the convergence rates and computational efficiency of solving partial differential equations (PDEs), making them highly effective for a wide range of complex problems. Recently, the combination of DDMs with neural networks for solving partial differential equations has attracted significant attention, as it can offer significant computational speed-ups robustness and enhanced accuracy for various deep learning models. These hybrid approaches are often classified based on the machine learning algorithms they incorporate \cite{heinlein2021combining,klawonn2024machine}. In this section, we aim to provide a  brief overview of combining DDMs with two prominent machine learning models.
	
	\subsection{Physics-informed neural networks}
	
	{In 1994, Dissanayake and  Phan-Thien \cite{dissanayake1994neural} introduced a neural network-based method that transforms solving PDEs into an unconstrained minimization problem, offering an efficient solution method. Later, Lagaris et al. \cite{lagaris1998artificial} developed a neural network approach for initial and boundary value problems by designing a trial solution that automatically satisfies boundary conditions and optimizes the network to solve the differential equation.}
	
	Based on those pioneering works, Physics-informed neural networks (PINNs) \cite{raissi2019physics} have recently been widely used to solve a variety of PDEs by using neural networks. Consider the following generic boundary value problem (BVP) on the domain $\Omega \subset \R^d$,
	\begin{equation} \label{eq1}
		\begin{aligned}
			\mathcal{A} [u] \;&=\; f \qquad {\rm in}\; \Omega, \\
			\mathcal{B} [u] \;&=\; g \qquad{\rm on}\; \partial \Omega, 
		\end{aligned}
	\end{equation}
	where $u: \Omega \to \R$ represents the solution, $\mathcal{A}$ is a differential operators, $\mathcal{B}$ is the boundary conditions operator, and  $f: \Omega \to \R$ is a given right-hand side function. The regularity of $u$ and $f$ depends on $\mathcal{A}$ and $\mathcal{B}$.
	
	To solve \eqref{eq1} using PINNs, we aim to find a neural network $u_p(x,\theta): \Omega\to \R$ that approximates the solution $u$. In the case of a fully connected feedforward network, $u_p(x,\theta)$ is given by
	\begin{equation*}
		u_p(x,\theta) = f_n \circ f_{n-1} \circ  \ldots \circ f_i \circ \ldots \circ f_1(x,\theta),
	\end{equation*}
	where $f_i(x,\theta_i) = \sigma_i(A_ix+b_i)$ for $i  =1,\ldots,n-1$ and $f_n(x,\theta_n) = A_nx+b_n$. Here, $A_i \in \mathbb{R}^{d_i \times d_{i-1}}$, $b_i \in \mathbb{R}^{d_i}$, $d_i$ is the number of neurons in the $i$-th layer and $\sigma_i$ is the activation function, $d_0$ and $d_n$ represent the input and output dimensions, respectively, which must align with the dimensions of the problem. Moreover, $\theta_i = (A_i,b_i)$ is the set of learnable parameters in the $i$-th layer, and $\theta = (\theta_1,\ldots,\theta_n)$ are all learnable parameters. 
	
	By constructing loss functions based on the residual terms of the equations in \eqref{eq1}, we employ a gradient descent method, such as the  Adam optimizer, to minimize the loss function. The loss function is a discretized version of $\mathcal{L}(\theta)$,
	\begin{equation}\label{pinn}
		\mathcal{L}(\theta) = \int_{\Omega} \big(\mathcal{A}[u_p](x,\theta) - f(x) \big)^2 dx + \lambda \int_{\partial \Omega} \big(\mathcal{B} [u_p](s,\theta) - g(s)\big)^2 ds,
	\end{equation}
	where $\lambda>0$ is a scalar weight used to balance the terms in the loss function.

	Alternatively, the neural network can take the form of the solution ansatz \cite{lagaris1998artificial,mcfall2009artificial,leake2020deep,lyu2020enforcing}. Take Dirichlet boundary conditions as an example. One can adapt the constraining operator	$\mathcal{C}$ such that 
	\begin{equation} \label{operator}
	\mathcal{C}	u_p(x,\theta) = L(x) u_p(x,\theta) + G(x) 
	\end{equation}
	satisfies the boundary conditions automatically. Here, $L(x)$ and $G(x)$ are some constructed functions chosen such that \eqref{operator}  meets the boundary conditions defined in \eqref{eq1}. For other types of boundary conditions, similar constructions can be found in \cite{leake2020deep,lyu2020enforcing}. Then, $\mathcal{L}(\theta)$ becomes
	\begin{equation}
		\mathcal{L}(\theta) =   \int_{\Omega} \big(\mathcal{A}[\mathcal{C}	u_p](x,\theta) - f(x) \big)^2 dx,
	\end{equation}
	which eliminates the need for penalty parameter $\lambda$ in \eqref{pinn}. However, now the operator $	\mathcal{C}$ has to be chosen appropriately.
	
	Extending the idea of PINNs, conservative PINNs \cite{jagtap2020conservative} adopt a non-overlapping domain decomposition method specifically designed for conservation laws by enforcing the flux continuity in the strong form along the subdomain interfaces. This approach is further developed in \cite{jagtap2020extended} as extended PINNs (XPINNs) to handle general PDEs and arbitrary space-time domains. Furthermore, $hp$-variational PINNs ($hp$-VPINNs) \cite{kharazmi2021hp} utilize piecewise polynomials as test functions within each subdomain, thereby improving the precision of solutions across the computational domain. Another approach, the deep domain decomposition method (DeepDDM) \cite{li2020deep}, leverages a classical fixed-point Schwarz iteration \cite{schwarz1870ueber} to decompose the boundary value problem. Finite basis PINNs (FBPINNs) were introduced in \cite{moseley2023finite}, where overlapping Schwarz domain decomposition is used to construct a global solution from localized NN-based functions.. This concept was further explored in \cite{dolean2022finite}, which applied additive, multiplicative, and hybrid iteration techniques to train FBPINNs. To solve problems with high frequency and multi-scale solutions, multilevel FBPINNs were introduced in \cite{dolean2024multilevel}, incorporating multiple levels of domain decomposition. Additionally, \cite{heinlein2024multifidelity} adopted a similar approach and enhanced the efficiency of time-dependent problems, and \cite{howard2024finite} developed a DDM for Kolmogorov-Arnold Networks (KANs), enabling the efficient training of KANs to solve multiscale problems.

	\subsection{Randomized Neural Networks}
	
	A major advantage of methods based on PINNs is that they provide a mesh-free algorithm, as they do not require discretization of the PDEs, and the differential operators in the governing PDEs are approximated using automatic differentiation. However, limitations in accuracy and training cost remain due to the non-convex optimization problem \cite{anderson2024elm}, (which can lead to local minima \cite{jagtap2020conservative}.)	
	
	For this reason, Pao et al. \cite{pao1992functional} introduced random vector functional-links (RVFLs), a method for single-hidden layer networks where connections to the hidden layer are fixed after random initialization, leaving only the output layer weights adjustable. This approach reduces training costs, improves computational efficiency \cite{pao1994learning}, and maintains universal approximation capabilities \cite{igelnik1995stochastic}. Extreme Learning Machines (ELMs) \cite{huang2006extreme} further advanced this concept by analytically computing output weights through solving a linear system, achieving fast learning speeds while maintaining universal approximation properties \cite{liu2014extreme}. A detailed survey of randomized neural networks (RaNNs) was later presented in \cite{gallicchio2020deep}.
	
	Taking the single-hidden layer as an example, a RaNN $u_r(x,W): \Omega\to \R$ is given by:
	\begin{equation} \label{rann}
		u_r(x,W) = W \cdot \sigma(Rx+b), 
	\end{equation}
	where  $R \in \R^{m\times d}$ and $b\in\R^m$ are randomly generated from a uniform distribution and fixed thereafter, and $W \in \R^{m}$ is the output weight that needs to be determined. The function $\sigma$ is the activation function.
		
	To solve  \eqref{eq1}  using \eqref{rann},  given a set of collocation points $\{x_i\}_{i=1}^{N} \subset \Omega $, and a constraining operator $\C$ defined in \eqref{operator},  the goal is  to find $W$ such that  ideally
	{\begin{equation} \label{elm}
	\mathcal{A} [L(x_i)W \cdot \big( \sigma(Rx_i+b)\big)+ G(x_i)]   = f(x_i)  \quad \forall i\in \{1, \ldots, N\}. 
	\end{equation}}If $\mathcal{A} $ is a linear operator, then we define matrix $E\in \R^{N\times m}$ with $E_{ik} =\mathcal{A} [L(x_i) \sigma(R_k x_i +b_k)]$ and a vector $T\in \R^N$ with $T_i = f(x_i)- \mathcal{A}[G(x_i)]$, then the above equation \eqref{elm} can be written as
	\begin{equation} \label{ls1}
		EW = T. 
	\end{equation}
	Here, $E$ is generally not invertible. Instead of using gradient-based algorithms like the Adam optimizer, the vector $W$ can be determined by the least-squares solutions of the system \eqref{ls1}. This approach provides a fast alternative to traditional neural network training techniques \cite{huang2006extreme}.
	
	The integration of local ELMs with DDMs was introduced in \cite{dong2021local}, offering improvements in both accuracy and efficiency. For problems with more complex geometries, a related approach called the random feature method, which utilizes random feature functions combined with overlapping domain decomposition, was developed in \cite{chen2022bridging}. A combination of RaNNs and the discontinuous Galerkin method was presented in \cite{sun2024local,sun2024local_W}, while the hybrid discontinuous Petrov-Galerkin method was combined with RaNNs in \cite{Dang2024local}.
	
	\subsection{Singular value decomposition }
	The use of the singular value decomposition (SVD) in neural networks has been explored to reduce the model size while keeping the accuracy in recent research. One approach is to examine the SVD of the weight matrices in multi-layer neural networks. In \cite{weigend1991effective}, the SVD of the weight matrix was analyzed to show how the dimension of the weight space changes during backpropagation training. A method for reducing redundant hidden units by analyzing the matrix rank of the three-layered feedforward neural network using SVD was studied in \cite{tamura1993determination}, enabling network simplification without increasing the training error. Then, an online algorithm for determining the smallest possible number of neurons based on an SVD was proposed in \cite{hayashi1993fast}. In \cite{xue2013restructuring}, for multi-layer neural networks with a large number of parameters, an SVD is applied to the weight matrices 
	$$\A_{b\times c} = U_{b\times e} \Sigma_{e\times e}V^T_{e\times c},$$ decomposing them into two smaller matrices, $U$ and $\Sigma V^T$ thereby reducing the number of parameters from $bc$ to $(b+c)e$ when $e$ is generally much smaller than $c$.
	
	On the other hand, once the hidden weights of a neural network are fixed \cite{psichogios1994svd}, a least-squares problem can be formulated to calculate the weights for the output layer, and then the SVD is computed on the hidden layer output matrix \cite{psichogios1994svd, abid2002new}. Psichogios and Ungar \cite{psichogios1994svd} proposed the SVD-NET algorithm to identify and eliminate redundant hidden nodes in a single hidden layer network. In \cite{abid2002new}, the authors find that the number of small singular values changes for different initialization of the weights. To address this, they analyzed the singular values across various initial weight settings and selected the optimal range of the initial weights. Teoh et al. \cite{teoh2006estimating} quantified the significance of the number of hidden neurons using SVD and proposed several possible methods for selecting the threshold parameter to eliminate the small singular values, illustrating that retraining the pruned network  by an appropriate threshold would be more efficient than training the original network, as the training process aims to maximize the linear independence of the reduced set of hidden layer neurons. 
	
	In this paper, we use RaNNs with domain decomposition to solve PDEs. Overlapping  Schwarz DDMs are utilized both to construct the global solution using RaNNs and to precondition the least-squares problem. In detail, the computational domain is partitioned into several subdomains, with RaNNs generated for each subdomain and initialized using a uniform distribution. Local window functions and the constraining operator are employed to construct a global approximate solution without the need for additional penalty terms.  Then, a linear system $H W = F$  is constructed by selecting a set of collocation points within the domain to enforce the PDEs at these points. To reduce the condition number of the linear system, cond$(H) = \frac{\sigma_{max}}{\sigma_{min}}$, where $\sigma_{max}$ and $\sigma_{min}$ represent the maximum and minimum singular values of $H$, respectively, we aim to eliminate small singular values, thereby increasing $\sigma_{min}$ to improve the conditioning of the system. The maximum singular value is bounded in terms of a coloring constant, see, e.g., \cite{widlund1988iterative}. Principal Component Analysis (PCA) \cite{abdi2010principal} is a statistical method used for dimensionality reduction and data compression by identifying principal components, which are linear combinations of the original variables, and eliminating correlations between the data.  Here, PCA is applied to the hidden layer output matrix of RaNNs using SVD, resulting in neural networks that require fewer parameters and significantly eliminate small singular values. Finally, we design overlapping Schwarz preconditioners to efficiently solve a least-squares problem using a preconditioned iterative method like Conjugate Gradient (CG) and Generalized Minimal Residual (GMRES) method.

	The paper is organized as follows: Section 2 introduces the fundamental concepts of  RaNNs integrated with overlapping Schwarz DDMs. Section 3 provides a detailed explanation of  additive and restricted additive Schwarz preconditioners. Section 4 illustrates how the PCA process is formulated for RaNNs. Section 5 presents numerical results to demonstrate the accuracy and efficiency of the proposed framework. The final section concludes the paper with a discussion of the findings and potential future work.
	
	\section{RaNNs with  domain decomposition}
	
	In overlapping Schwarz domain decomposition methods, the computational domain $\Omega$ is divided into $J$ overlapping subdomains $\{\Omega_j\}_{j=1}^J$ classically. The global discretization space is often defined first, and then the local spaces are constructed as subspaces according to the domain decomposition \cite{smithdomain,toselli2006domain}. Here, we define local spaces $V_j$  of neural network functions  and construct the global discretization space $V$ from these local spaces using window functions. 
	
	Let $u_j(x,W_j)$ be a RaNN with a fully connected feedforward structure defined on subdomain $\Omega_j$ given by
	\begin{equation} \label{net}
		u_j(x,W_j) = W_j \cdot f_j^{n-1} \circ  \ldots \circ f_j^i \circ \ldots \circ f_j^1(x) := W_j \cdot \Phi_j(x),
	\end{equation}
	where $\Phi_j(x) =  f_j^{n-1} \circ  \ldots \circ f_j^i \circ \ldots \circ f_j^1(x)$ is fixed and contains no learnable parameters. Each layer is defined as $f_j^i(x) = \sigma_i(R_ix+b_i)$, for $i  =1,\ldots,n-1$, and $R_i, b_i$ are randomly generated from a given uniform distribution and remain fixed.  Therefore, only $W_j \in \R^m$ needs to be determined.
	
	Then, we can define a function space of RaNNs as follows:
	\begin{equation} \label{space_vj}
		V_j = \text{span}\{\phi_j^1(x),\cdots, \phi_j^m(x)\},
	\end{equation}
	where $\phi_j^k(x)$ denotes the output of the $k$-th neuron in the $n-1$ layer, and they are the components of $\Phi_j(x)$, that is, $\Phi_j = (\phi_j^1,\ldots,\phi_j^m)$. $\{\phi_j^k\}_{k=1}^{m}$ can be regarded as a set of basis functions. These functions generally have global support, hence, we need to localize them using window functions.
	
	Therefore, we introduce  $\{\omega_j\}_{j=1}^{J}$ as smooth window functions that form a partition of unity, confining the neural networks to their respective subdomains. In particular, we choose $\omega_j$ such that 
	\begin{equation}
		\mathrm{supp}(\omega_j)  \subset \overline{\Omega}_j  \quad \text{and} \quad \sum_{j=1}^{J}\omega_j \equiv 1  \quad \text{on} \;\Omega.
	\end{equation}

\begin{figure}[!ht] 		
	\centering
	\includegraphics[scale=0.5]{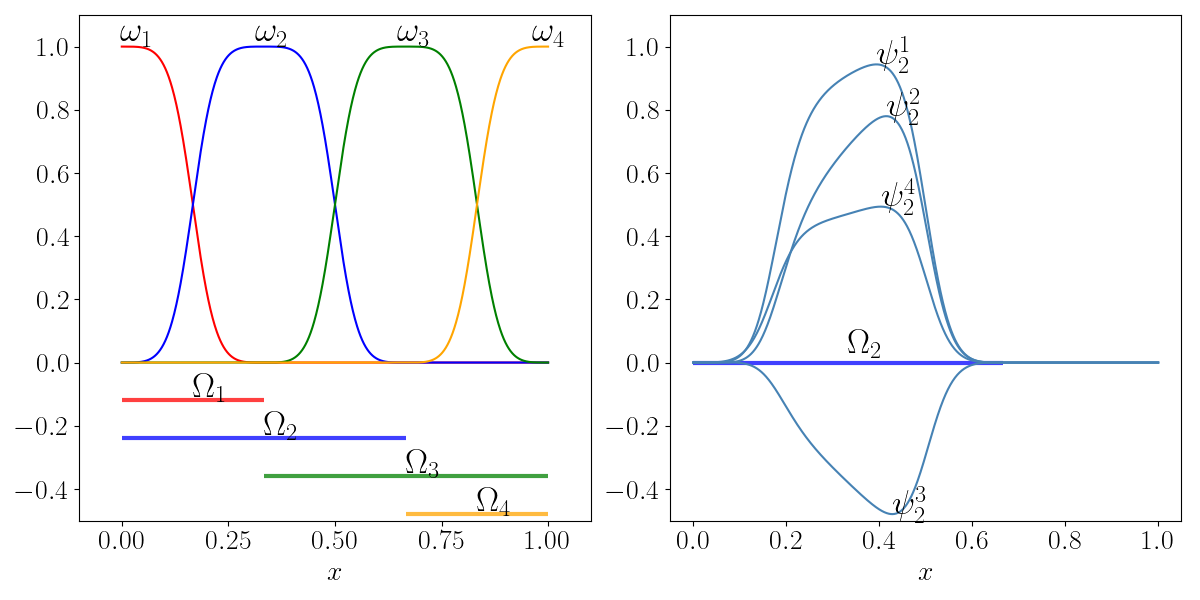}
	\caption{Example of using RaNNs with overlapping Schwarz DDM to solve PDEs  in one dimension. The domain is divided into four subdomains, and RaNNs with one hidden layer and four hidden nodes are used for each subdomain.  The window functions $w_j$, each with local support within $\Omega_i$, are displayed in the left figure. Exemplary local basis functions  $\psi_j^k(x)$ for the second subdomain are shown in the right figure.}
	\label{domain}
\end{figure}

	Similar to multilevel FBPINNs paper \cite{dolean2024multilevel}, we can define the global space of RaNN functions as
	\begin{equation} \label{space_v}
		V = \sum_{j = 1}^{J} \tilde{V}_j.
	\end{equation}
	Here, the localized neural network function space is defined as
	\begin{equation} \label{space_local}
		\tilde{V}_j =\text{span}\{\psi_j^1(x),\cdots, \psi_j^m(x)\},  
	\end{equation}
	where
	\begin{equation}
		\psi_j^k(x) = \omega_j \phi_j^k(x), \quad j = 1,\ldots,J, \; k = 1,\ldots,m.
	\end{equation}
	can be regarded as local basis functions. In other words, $\tilde{V}_j = \omega_j V_j$. An example of window functions $\omega_j$ and local basis functions $\psi_j^k$ is shown in Figure \ref{domain}.
	
	Thus, the approximated solution $\hat{u}(x)$ in $V$ can be represented as a combination of contributions from networks on all subdomains, multiplied by a constraining operator $\C$ defined in \eqref{operator} to enforce the boundary conditions,
	\begin{equation} \label{netdd}
		\hat{u}(x) =  \C\sum_{j=1}^{J} \omega_j \,	u_j(x) = L(x) \sum_{j=1}^{J}\omega_j W_j \cdot \Phi_j(x) + G(x) =  \sum_{j=1}^{J}W_j \cdot \big(L(x)\omega_j \Phi_j(x)\big) + G(x) .
	\end{equation}	
	where $L(x)$ and $G(x)$ are constructed functions related to the domain $\Omega$, designed such that equation \eqref{netdd} satisfies the boundary conditions defined in \eqref{eq1}.

	By substituting equation \eqref{netdd} into \eqref{eq1}, and selecting a set of collocation points $\{x_i\}_{i=1}^{N} \subset \Omega $, we aim for the approximated solution  $\hat{u}(x)$ to satisfy the PDEs within the global domain.  We need to find $W_j$, $j = 1,\ldots,J$, such that 
	\begin{equation}
	 \mathcal{A} \Big[\sum_{j=1}^{J}W_j \cdot \big(L(x_i)\omega_j \Phi_j(x_i)\big) \Big]  = f(x_i) - \mathcal{A}[G(x_i)],  \quad \forall i\in \{1, \ldots, N\}. 
	\end{equation}
	When $\mathcal{A}$ is a linear operator, we have
	\begin{equation} \label{sum}
	\sum_{j=1}^{J} 	W_j \cdot \mathcal{A} \big[L(x_i) \omega_j  \Phi_j(x_i) \big]  = f(x_i) - \mathcal{A}[G(x_i)],  \quad \forall i\in \{1, \ldots, N\}. 
	\end{equation}
	The above equation yields
	\begin{equation}\label{system}
	[H_1,H_2,\ldots,H_J](W_1,W_2,\ldots,W_J)^T = F,
	\end{equation}
	where $H_j\in \R^{N\times m}$ with $(H_{j})_{i,k} = \mathcal{A} [L\omega_j  \phi_j^k(x_i)]$, and $F\in \R^N$ with $f(x_i) - \mathcal{A}[G(x_i)]$. 
	
	\begin{remark}
		In the case where $\mathcal{A}$ is a nonlinear operator, nonlinear iterative methods such as Picard iteration and Newton's method can be employed, as discussed in \cite{shang2023Randomized}. In each iteration step, a linear least-squares problem is solved.
	\end{remark}
	For simplicity, we rewrite the linear system \eqref{system} as
	\begin{equation} \label{sys}
		H W = F,
	\end{equation}
	where $H = [H_1,H_2,\ldots,H_J]$ and $W = (W_1^T,W_2^T,\ldots,W_J^T)^T$.
	In general, the matrix $H$ is not square (and not symmetric), and we solve it using the least-squares method. To solve the least-squares problem
		\begin{equation*}
			\min \Vert H W - F\Vert_2^2,
		\end{equation*}
		the first-order optimality condition leads to the normal equation:
		\begin{equation} 
			H^T H W = H^T F.
	\end{equation}Here, $H^T H \in \R^{p \times p}$ is a sparse  symmetric matrix, and $p=Jm$. The condition number of $H^TH$ is is given by cond($H^TH$) = cond($H$)$^2$. 

	A large condition number arises from the linear dependencies among the local basis functions $\psi_j^k$, which can result in a rank-deficient matrix. This, in turn, leads to small singular values and an increased squared condition number for the normal equations, if used. We propose additive and restricted additive Schwarz preconditioners, and use Conjugate Gradient (CG) and Generalized Minimal Residual (GMRES) methods to solve the problem, respectively.

	\section{Overlapping Schwarz preconditioners}
	
	{Consider an example using RaNNs with overlapping Schwarz DDM  to solve the Laplace equation \eqref{exlapla} in two dimensions. The domain is divided into 16 subdomains, and RaNNs with one hidden layer and $m$ hidden nodes are employed on each subdomain.  Figure \ref{eigens} displays the sparsity patterns of the matrix for different values of $m$ when using $40 \times 40$ collocation points, showing that, as more parameters are introduced to improve accuracy (i.e., as $m$ increases), the condition number of the matrix $H$ and $H^T H$ increases as well. The distribution of eigenvalues of the matrix $H^T H$ is shown in Figure \ref{eigens}, with the condition number being cond($H^TH$) = $10^{16}$. We can observe that the minimum eigenvalue is around $10^{-10}$ and the maxmum eigenvalue is around $10^{6}$, which would lead to slow convergence for iterative methods like CG and GMRES \cite{luenberger1973introduction}. Therefore, an effective preconditioner is needed to improve the convergence rate and address the ill-conditioning of the system.}
	
	\begin{figure}[!ht]
		\begin{center}
			\subfigure[ $m  = 16$, cond($H^TH$) = $10^{12}$.]{
				\centering
				\includegraphics[width=3.0in]{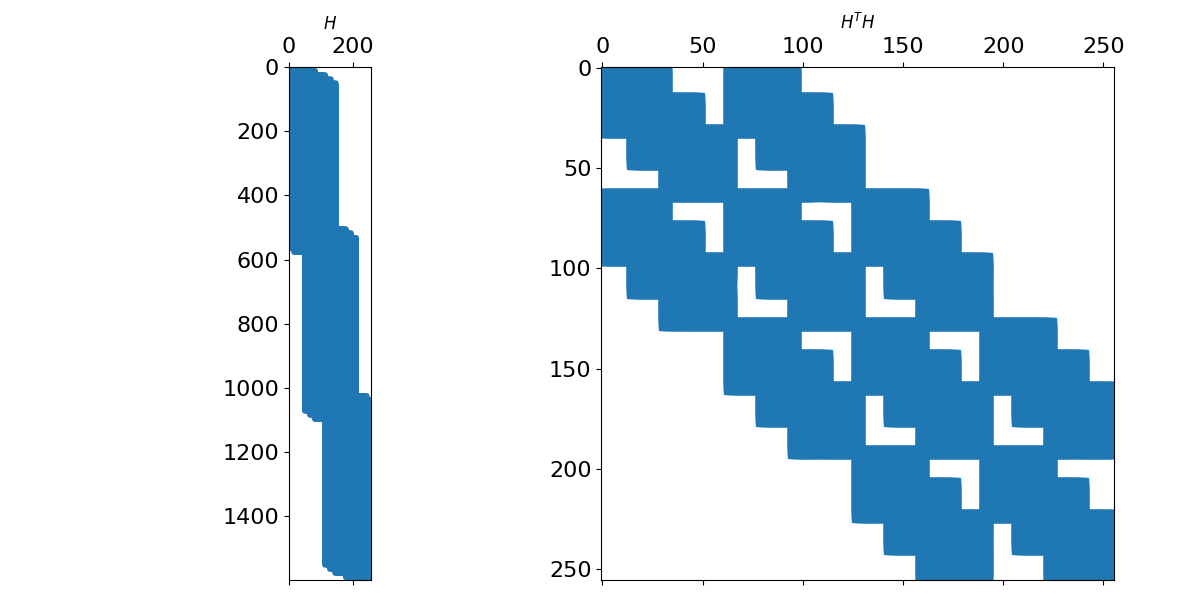}
			}
			\subfigure[$m  = 32$, cond($H^TH$) = $10^{16}$.]{
				\centering
				\includegraphics[width=3.0in]{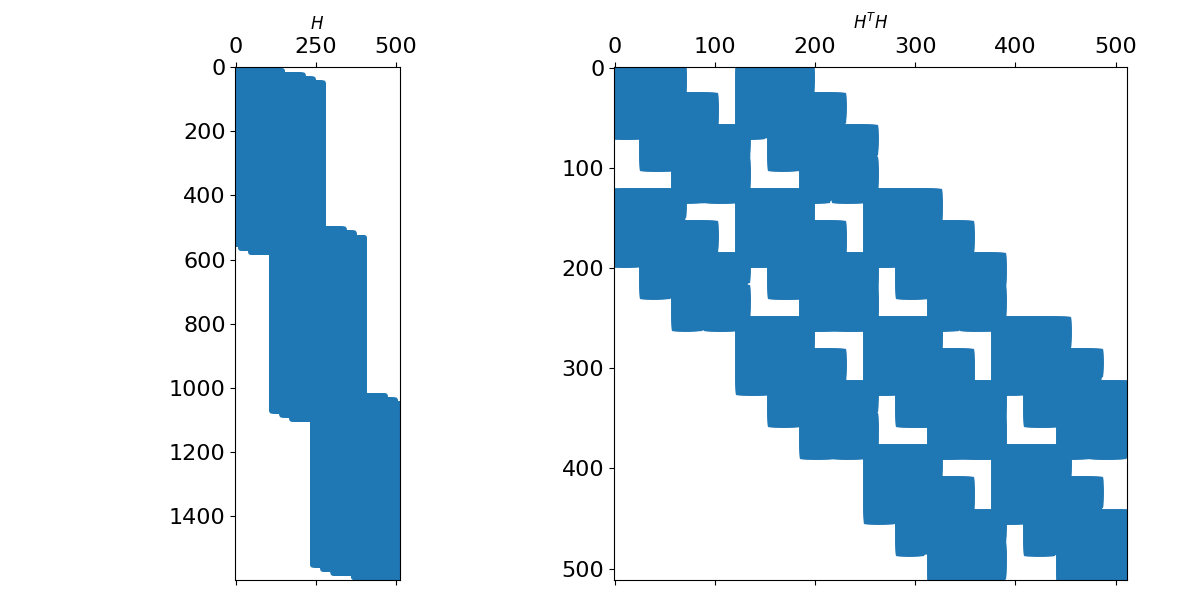}
			}
		\end{center}
		\vspace*{-15pt}
		\caption{Different values of $m$ are tested in (a) and (b), where the condition numbers are shown for $m = 16$ and $m =32$, respectively.}
		\label{matrix}
	\end{figure}

	\begin{figure}[!htbp]
		\includegraphics[width=6.2in]{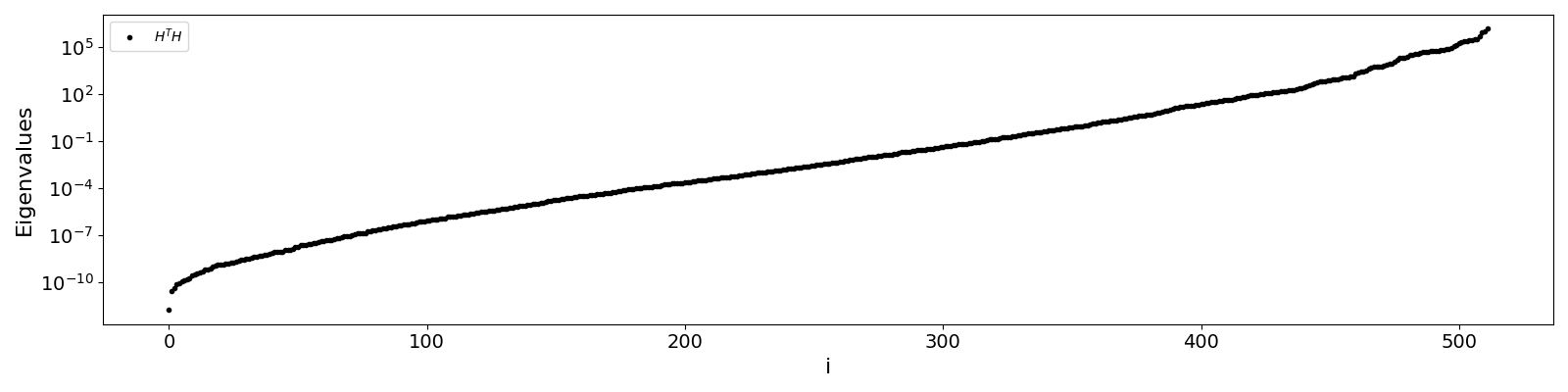}	
		\caption{The distribution of all 512 eigenvalues of the matrix $H^T H$ is shown.}
		\label{eigens}
	\end{figure}

	As Chen et al. mentioned in \cite{chen2024optimization}, solvers for least-squares problems have not been well-developed, particularly for ill-conditioned problems with sparse structures. The condition number of the matrix in RaNN-based methods is extremely large \cite{shang2023Randomized}. Therefore, we propose using preconditioners to accelerate convergence in iterative methods such as the CG and GMRES method. Schwarz preconditioners are widely employed in DDMs to enhance the efficiency of iterative solvers for PDEs; (see \cite{chan1994domain,toselli2006domain} for instance.) Therefore, we focus on constructing overlapping Schwarz preconditioners specifically for RaNN-based methods.

	Based on \eqref{space_v} and \eqref{space_local},we can define restriction operators
	\begin{equation*}
		R_i : V\rightarrow \tilde{V_i},  \quad R_i(\sum_{j=1}^{J} \omega_j u_j) = \omega_i u_i, \quad i = 1,\ldots,J,
	\end{equation*}
	which take a function, expressed as a sum of functions from $\tilde{V}_i$, and maps it onto the term in $\tilde{V}_i$. 
	Similarly,  extension operators are
	\begin{equation*}
		R_i^T : \tilde{V_i} \rightarrow V,  \quad R_i^T(\omega_i u_i) = \omega_i u_i, \quad i = 1,\ldots,J.
	\end{equation*}
	
	Then, an additive Schwarz (AS) preconditioner for $H^T H$ is given by
	\begin{equation}
		M^{-1}_{AS} = \sum_{i = 1}^{J} R_i^T A_{i}^{-1} R_i,  
	\end{equation}
	where  $A_{i} = R_i H^TH R_i^T$, for $i = 1,\ldots,J$. 	In practice, let $R_i$ be a $q \times p$ Boolean matrix, and let $s_j$ denote the indices of the local basis functions $\psi_j^k(x)$ that lie within the interior of $\Omega_j$. We can construct the matrix such that $R_i[:, s_j] = I_q$, where $ I_{q} \in \R^{q\times q}$ is the identity matrix, and $R_i^T$ is constructed analogously.
	
	By reducing the communication between different subdomains, a more effective preconditioner called the restricted additive Schwarz (RAS) preconditioner was proposed in \cite{cai1999restricted}. After introducing the algebraic partition of unity matrix decribed, e.g., in \cite{dolean2015introduction}, denoted by $D_i \in \R^{n_i\times n_i}$, for $i = 1,\ldots, J$, 
	the RAS preconditioner can be represented as:
	\begin{equation}
		M^{-1}_{RAS} = \sum_{i = 1}^{N} {R}_i^T D_i A_{i}^{-1} R_i, 
	\end{equation}
	Here, $D_i$ is a diagonal nonnegative matrix that satisfies the relation
	\begin{equation*}
		\sum_{i=1}^{J} R_i^T D_i R_i = I_{p}.
	\end{equation*}  
	where $ I_{p} \in \R^{p\times p}$ is the identity matrix.
	
	The inverse $A_{i}^{-1}$ is computed using the QR decomposition, which is numerically stable,
	\begin{equation*}
		A_i = Q_iP_i,   \quad i = 1,\ldots,J,
	\end{equation*}
	where $Q_i$ is an orthogonal matrix  with  $Q_i^T = Q_i^{-1}$ and  $P_i$ is an upper triangular matrix. Here, we call it $P_i$ as we already use the symbol $R_i$. Then, we have
	\begin{equation*}
		A_{i}^{-1} = P_i^{-1} Q_i^T, \quad i = 1,\ldots,J,
	\end{equation*}	
	where $P_i^{-1}$ denotes the inverse of $P_i$ when it is invertible; otherwise, the pseudo-inverse $P^{+}$ is used, with both computed using Gaussian elimination.
	
	\begin{remark}	
		To (naturally) obtain a matrix $H$ with a relatively low condition number, we can select parameters that help avoid linear dependencies of $\psi_j^k(x)$, as linear dependencies tend to introduce very small eigenvalues, hence increasing the condition number. One effective approach is to decrease the similarity of data representations of $\psi_j^k(x)$ in the overlapping regions of subdomains, This can be accomplished by ensuring there are more data points $N$ than hidden neurons $m$. 
		
		Additionally, combining RaNNs with overlapping Schwarz DDMs results in a sparse matrix as shown in Figure \ref{matrix}, which takes the form of a block-sparse matrix. This structure allows us to construct the overlapping Schwarz preconditioners in a natural way. Without DDMs, the matrix would be full, making it challenging to design an effective preconditioner.

	\end{remark}
	\section{PCA process }
	
\begin{figure}[!htbp]
		\begin{center}
			\subfigure[Original RaNN]{
				\centering
				\includegraphics[scale=0.43]{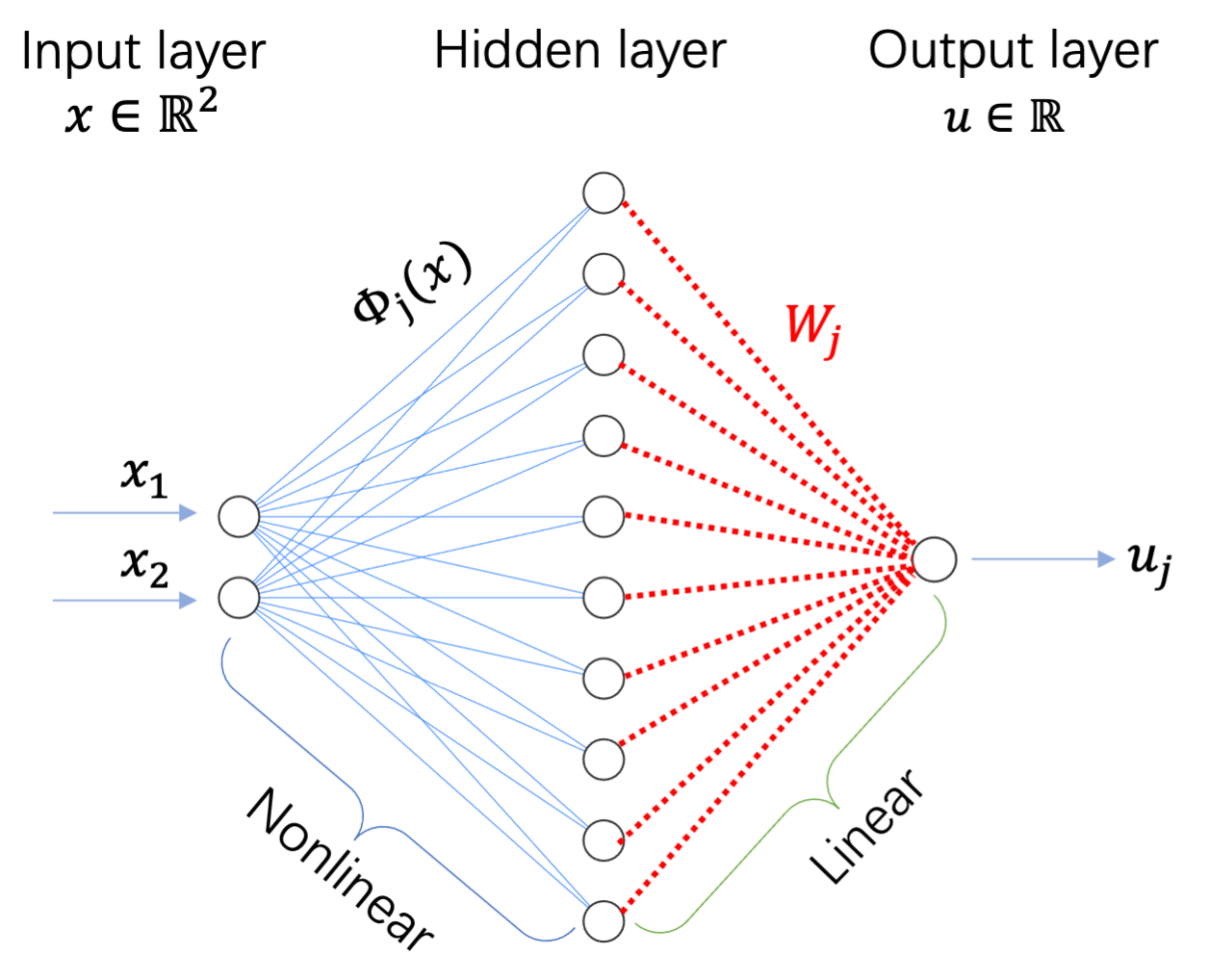}
			}
			\hspace{7mm}
			\subfigure[RaNN with a PCA process]{
				\centering
				\includegraphics[scale=0.43]{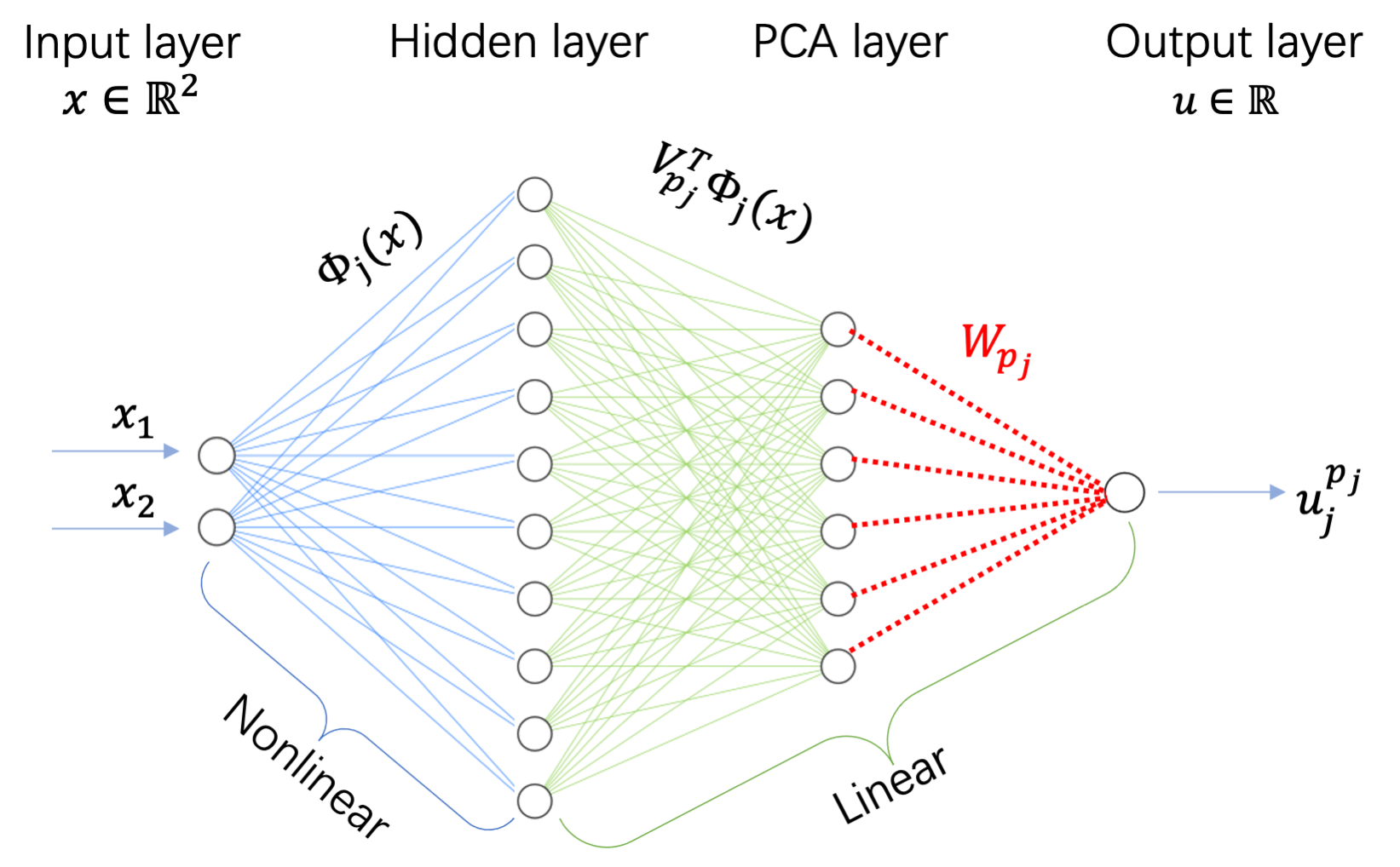}
			}
		\end{center}
		\vspace*{-15pt}
		\caption{Network structure of the original RaNN $u_j:\mathbb{R}^{2}\rightarrow \mathbb{R}$ and new output $u_j^{p_j}$ after a PCA process, the solid blue line represents the parameters of the neural network that are randomly initialized and fixed thereafter, and the solid green line represents the PCA process, which involves a matrix multiplication with $V_{p_j} ^T$. The dotted red line indicates the parameters that need to be determined.}
		\label{network}
	\end{figure}

	{The phenomenon known as condensation of neural networks is discussed in \cite{luo2021phase,zhou2022towards}, where it is observed that for networks with random initialization, the weights from an input node to hidden neurons quickly converge to similar values after training. This results in multiple hidden neurons can be replaced by an effective neuron with low complexity \cite{zhou2022towards}. 
	
	For RaNNs, we apply a PCA process as shown in Figure \ref{network} (b), which preserves the principal components of the neural network, and reduces the number of solvable parameters while maintaining much of the network's approximation capability.  This approach helps to mitigate linear dependencies within the system, resulting in a lower condition number.}
	
	 From the construction  of $H$, we know 
	\begin{equation}
		H = 	[H_1,H_2,\ldots,H_J], \quad (H_j)_{i,k} = \mathcal{A} [L\omega_j  \phi_j^k(x_i)].
	\end{equation}
	Here, for a given problem \eqref{eq1} with operator $\mathcal{A}$, we can either perform an SVD on $H_{j}$ directly or analyze the output matrix of the local basis functions $\psi_j^k(x)$ denoted by
	\begin{equation}
		(\Psi_j)_{i,k} =\psi_j^k(x_i),   \; j = 1,\ldots,J.
	\end{equation}
	The latter approach enables us to focus on the key components in the network output that contribute to the solution approximation. It is problem independent and can be done in an offline process, ensuring a more efficient solution for different problem settings after initialization of the weights.
	
	After performing a reduced singular value decomposition (SVD) on $\Psi_j$, we obtain \begin{equation*}
		\Psi_j = U_j \Sigma_j V_j^T, \quad j = 1,\ldots,J,
	\end{equation*}
	where $k = \min(N_1,m)$, $U_j \in \R^{N_1 \times k}$, $\Sigma_j \in \R^{k \times k}$,  and $V_j \in \R^{m \times k}$.
	
	Then, we simply  select a constant threshold parameter $\tau$ to identify the singular values greater than $\tau$, denoting their count as $p_j$ on subdomain $\Omega_j$. Next, we choose the first $p_j$-$th$ columns of  $V_j$ to form the truncated matrix $V_{p_j} \in \R^{m \times p_j}$, respectively. 
	
	PCA process can be performed based on the SVD of $\Psi_j$, and the truncated hidden layer output matrix $T_j \in \R^{N_1 \times p_j} $ can be obtained by:
	\begin{equation} \label{pca}
		T_j = \Psi_j V_{p_j} , \quad j = 1,\ldots,J.
	\end{equation} 

	Now, we integrate a PCA process into the original RaNN, and the new output of the RaNN is given by:
	\begin{align} \label{pca_net}
	u_j^{p_j}(x)=& \sum_{k=1}^{p_j} W_{p_j}^k \sum_{i=1}^{m}\Phi_j^i(x) (V_{p_j})_{i,k} \notag \\
	=&\sum_{k=1}^{p_j} W_{p_j}^k \sum_{i=1}^{m} (V_{p_j}^T)_{k,i} \Phi_j^i(x) \notag \\ 
	=& W_{p_j} \cdot [V_{p_j} ^T \Phi_j(x)]  .
	\end{align}Here, $p_j$ denotes the number of effective neurons, and  $W_{p_j} \in \R^{p_j}$ are the new weights that need to be determined, and $W_{p_j}^k$ denote the $k$-$th$ row of $W_{p_j}$. The total number of parameters to be determined is given by  $\sum\limits_{j=1}^{J} p_j$, which does not exceed $p=mJ$.
	
	In Figure \ref{network} (a), the network structure of the original RaNN is illustrated, and in Figure \ref{network} (b), the PCA process is shown, where $\Phi_j(x)$ represents  the output of the nonlinear part. The term $V_{p_j} ^T\Phi_j(x)$ is referred to as the PCA layer, which is linear and only involves matrix multiplication, and does not introduce any additional unknown parameters.

	Similar to \eqref{sum}, we obtain the following linear system, aiming  to find $W_{p_j}$, $j = 1,\ldots,J$, such that  
	\begin{equation}  \label{system_new}
		\sum_{j=1}^{J} 	W_{p_j} \cdot \mathcal{A} \big(L \omega_j  [V_{p_j} ^T\Phi_j(x)]  \big)  = f(x_i) - \mathcal{A}(G(x_i)),  \quad \forall i\in \{1, \ldots, N\}. 
	\end{equation}
	Finally, we can construct the AS and RAS preconditioners and use the CG or GMRES method to solve  \eqref{system_new}, respectively. The whole algorithm is summarized in Algorithm 1.	
	\begin{table}[!ht]
		\renewcommand{\arraystretch}{1.2}
		\setlength\tabcolsep{0mm}
		\begin{tabular}{l}
			\toprule
			{\bf  Algorithm 1} Overlapping Schwarz preconditioners for RaNNs with domain decomposition. \\ \midrule
			\; Step 1. Use overlapping Schwarz DDM to divide the domain $\Omega$ into $J$ overlapping subdomains $\{\Omega_j\}_{j=1}^J$.  \\
			\; Step 2. Initialize $J$ RaNNs with network architecture $u_j :\Omega_j \rightarrow \mathbb{R}$ using a uniform distribution. \\
			\; Step 3. Construct a global approximate solution by applying window functions and a constraint operator. \\
			\; Step 4. Use collocation points $\{x_i\}_{i=1}^{N}$ and threshold parameter $\tau$ to perform PCA process.  \\
			\; Step 5. Use collocation points $\{x_i\}_{i=1}^{N}$ to assemble the linear system \eqref{system_new}. \\
			\; Step 6: Construct the AS and RAS preconditioners and solve the problem \eqref{system_new} with iterative method. \\
			\; Step 7. Solve the network parameters $W_{p_j}$ to  $u_j^{p_j}$, $ j = 1,\ldots,J.$\\ \bottomrule
		\end{tabular}
	\end{table}

	\begin{remark}
		A weight initialization strategy based on frequency information for constructing the initial RaNN to solve PDEs with varying solutions is presented in \cite{dang2024adaptive}. Additionally, \cite{zhang2024transferable} proposes a transferable neural network model based on shallow networks, making it easier to adapt the neural feature space across different PDEs in various domains. In this paper, since the inputs to each subdomain network are normalized, we use a uniform distribution for initialization. Further investigations into various weight initialization strategies will be explored.
		
		In practice, the PCA process is influenced by various factors, including the initialization of weights in neural networks, the number of collocation points and the choice of threshold parameter $\tau$. The impact of these factors is to reduce the condition number of the system, further investigation is needed to better understand their effects. This is out of the scope of the current work, and will be explored in future investigations.
	\end{remark}
	\section{Numerical examples}	
	All results are obtained using an NVIDIA GeForce RTX 3090 GPU. We use $scipy.linalg.lstsq$ in Python, relying on QR decomposition to directly solve the least-squares problem, and adopt the CG method and its variants, as well as the GMRES method from $scipy.sparse.linalg$, to solve the normal equation corresponding to the least-squares problem. We average the results over ten experiments for each problem to account for the randomness in the methods.

	\begin{example} \label{ex1}
		We consider the multi-scale Laplacian problem  from \cite{dolean2024multilevel}  with  Dirichlet boundary condition, 
		\begin{equation}
			\begin{aligned}\label{exlapla}
				-\Delta u &= f \quad \;\;\;\,{\rm in} \; \Omega = [0,1]^2,  \\
				u &= 0\qquad{\rm on} \; \partial \Omega.
			\end{aligned}
		\end{equation}
		and the exact solution  given by 
		\begin{equation*}
			u(x_1,x_2) = \frac{1}{n}\sum_{i = 1}^{n}\sin(\omega_i\pi x_1)\sin(\omega_i\pi x_2),
		\end{equation*}
		where $n$ = 1,2,$\cdots$6, and $\omega_i = 2^i, i = 1,2,\cdots n$, and $f =\frac{2}{n}\sum\limits_{i = 1}^{n}(\omega_i \pi^2)\sin(\omega_i\pi x_1)\sin(\omega_i\pi x_2)$.
	\end{example} 
	
	In this example, the domain $\Omega$ is divided into $J = l \times l$ overlapping subdomains, and a uniform rectangular domain decomposition is constructed.  Each dimension is defined as follows:
	\begin{equation*}
	\Omega_j = \left[\frac{(j-1)-\delta/2}{l-1}, \frac{(j-1)+\delta/2}{l-1} \right] \quad l\geq 1,
	\end{equation*}
	where $\delta > 1$ is defined as the ratio of the subdomain size and the overlap width between subdomains. 
	
	Same as in \cite{dolean2024multilevel}, the subdomain window functions are given by 
	\begin{align}
		\omega_j = \frac{\hat{\omega}_j }{\sum_{j=1}^{J^{(l)}}\hat{\omega}_j } \quad \text{where}   \quad \omega_j= \Pi_i^2[1+\cos(\pi(x_i-\mu_{i})/\sigma_{i})]^2,
	\end{align}
	where $\mu_{i} = (i -1)/(l -1)$ and $\sigma_{i} = (\delta/2)/(l-1)$ represent the center and half-width of each subdomain along each dimension, respectively.
	The constraining operator $\mathcal{C}$  is choosen using 
	\begin{equation} \label{cons1}
		L(x)= \tanh(x_1/\rho)\tanh((1-x_1)/\rho)\tanh(x_2/\rho)\tanh((1-x_2)/\rho), \quad G(x) = 0
	\end{equation}
	with $\rho = (1/2)^n$, ensuring that $\mathcal{C}u(x) = L(x) u(x)+ G(x) $ automatically satisfies the boundary condition. 
	
	First, we consider using a $4 \times 4 $ domain decomposition with $40 \times 40$ collocation points across the domain $\Omega$ to sove the problem \eqref{exlapla} for $n = 2$.  The overlap ratio $\delta$ is fixed by 2.	RaNNs with one hidden layer and 16 hidden nodes on each subdomain are initialized by a uniform distribution $\mathcal{U}(-1,1)$ for weights and bias parameters. Inputs to each subdomain network are normalized to the range $[-1,1]^2$ within the respective subdomains. The number of parameters that need to be solved is denoted as the degrees of freedom (DoF) and is given by DoF = $mJ$.  We first illustrate the performance of the overlapping Schwarz preconditioners without applying the PCA process. We construct the AS  preconditioner and denoted as $M^{-1}_{AS}$. For the RAS preconditioner,  we consider two choices for the coefficients of $D_i$ in $M^{-1}_{RAS}$: one is to use $(D_i)_{jj} = 1/\mathcal{M}_j$ denoted by $M^{-1}_{SAS}$ (scaled) \cite{dolean2015introduction}, where $\mathcal{M}_j$ is defined as the number of overlapping subdomains for each node $j$; the other is to use a Boolean matrix and denoted by $M^{-1}_{RAS}$ (restricted).
	
	In Table \ref{table1a}, $\kappa(H^TH)$ and $\kappa(M^{-1}H^TH)$ denote the condition number of the matrices $H^TH$ and $M^{-1}H^TH$ for different kinds of preconditioners, $\sigma_{min}$ and $\sigma_{max}$ denote the minimum and maximum singular values, respectively. and  $\lambda_{min}$ and $\lambda_{max}$ denote the minimum and maximum real parts of the eigenvalues, respectively. As for the iterative solvers, $M^{-1}_{AS}$ is symmetric while $M^{-1}_{SAS}$, $M^{-1}_{RAS}$ are not. Therefore, the GMRES method is applicable to all of them, as it does not require symmetry. The iteration is stopped when the relative residual norm of the solution falls below a fixed tolerance, i.e., \begin{equation*}
		\Vert M^{-1}H^THW - M^{-1}H^TF \Vert_{L^2} / \Vert M^{-1}H^T F \Vert_{L^2} \leq 10^{-5}.
	\end{equation*} Then, we report the number of iterations and relative $L^2$ error by iter and $e_{L^2}$, respectively. In particular, $e_{L^2}$ is computed approximately on $M = 350 \times 350$ uniformly-spaced test points by $ \Vert \hat{u}(x_i) - u(x_i) \Vert_{l^2} / \Vert {u}(x_i)\Vert_{l^2} $. From Table \ref{table1a}, we observe that the condition number of the matrix $H^TH$ is approximately $10^{12}$, and all three preconditioners significantly reduce the condition number around four orders of magnitude, and enable GMRES to converge in siginificantly fewer number of iteration steps compared to the case without preconditioning.
	
	\begin{table}[!ht]	
		\centering  
		\renewcommand{\arraystretch}{1.1}
		\setlength\tabcolsep{2.6mm}
		\caption{{Relative $L^2$ error and number of iterations for different preconditioners by GMRES in Example \ref{ex1}.}}	\label{table1a}	
		\begin{tabular}{ccccccccccc}  
			\toprule 
			DoF & $N$ &  $\kappa(H^TH)$ & $M^{-1}$  & $\kappa(M^{-1}H^TH)$ & $\sigma_{min}$ & $\sigma_{max}$ & $\lambda_{min}$ & $\lambda_{max}$  & iter  & $e_{L^2}$ \\ \midrule
			\multirow{4}{*}{256} & 	\multirow{4}{*}{1600} & \multirow{4}{*}{$10^{12}$ }  & 
			$None$ & $ 10^{12}$ & $10^{-7}$ & $10^{5}$ & $10^{-7}$ & $10^{5}$  & 2560 & 8.94e-2\\
			&&& $M^{-1}_{AS}$ & $10^{8}$ & $10^{-4}$ & $10^{4}$ & 4.6 & 16 & 13 &  6.44e-3 \\
			&&& $M^{-1}_{SAS}$ & $10^{8}$ & $10^{-5}$ & $10^{3}$   & 0.5 & 1.6 & 10 &  6.45e-3 \\
			&&& $M^{-1}_{RAS}$ & $10^{8}$ & $10^{-5}$ & $10^{3}$  & 0.1 &  2.1 & 31 &  6.41e-3\\
			\bottomrule
		\end{tabular}  
	\end{table}

	\begin{figure}[!htbp]
		\includegraphics[width=6.45in]{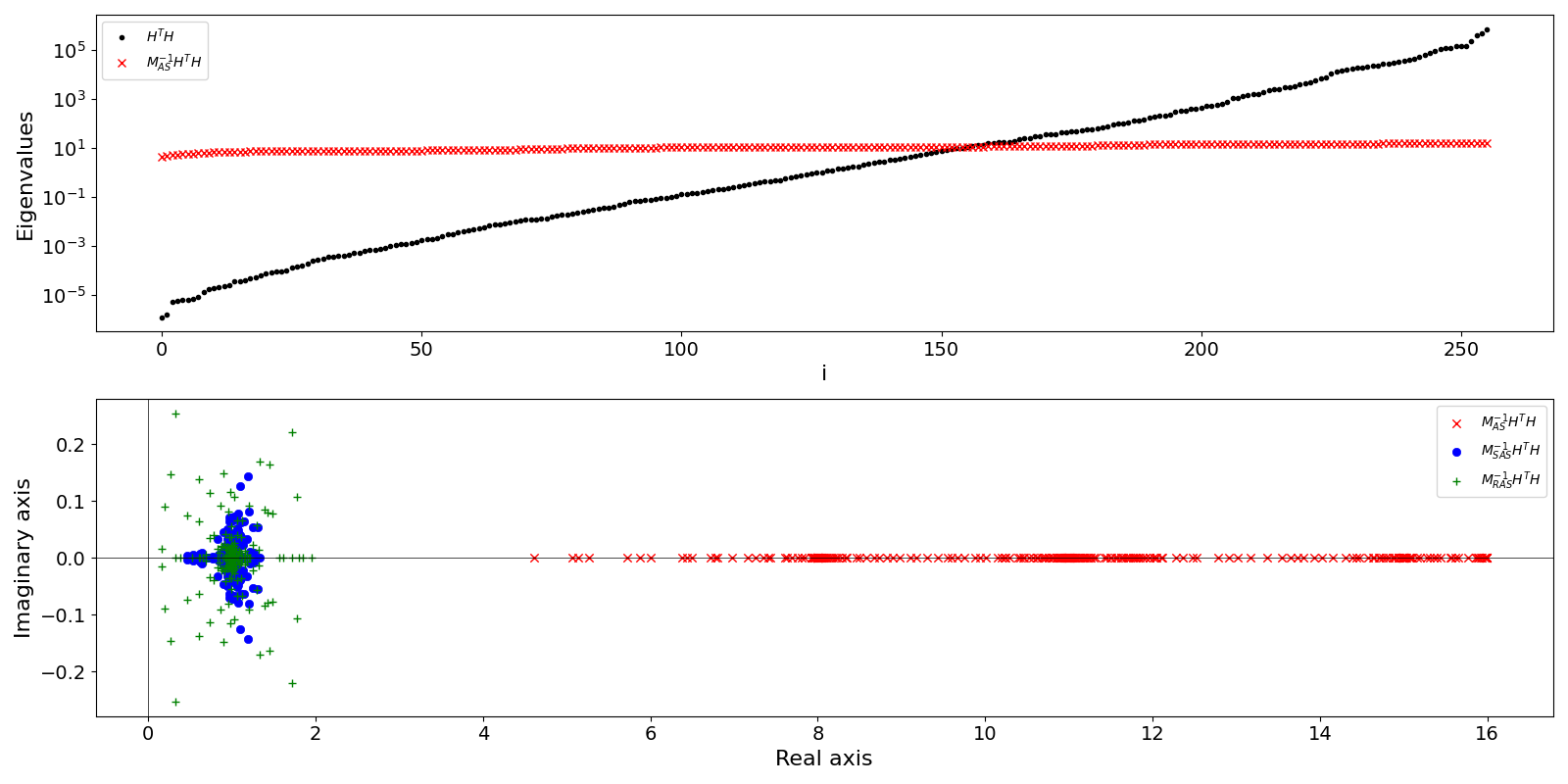}	
		\caption{The distribution of eigenvalues of the preconditioned matrix with different preconditioners  using $4 \times 4 $ domain decomposition in Example \ref{ex1}. The top figure compares the distribution of eigenvalues for the AS preconditioner against the distribution without preconditioning ordered in non-decreasing order, while the bottom figure illustrates the eigenvalue distributions for various preconditioners in the complex plane.}
		\label{figure1}
	\end{figure}

	The top of Figure \ref{figure1} shows the distribution of 256 eigenvalues for $H^T H$ and $M^{-1}_{AS} H^T H$, while the bottom shows the eigenvalue distribution in the complex plane for different preconditioners. We observe that the distribution of the eigenvalues of the AS preconditioned matrix becomes more concentrated, effectively shifting the minimum eigenvalue from $10^{-7}$ to $4.6$, with the maximum eigenvalue being reduced form $10^5$ to 16. In classical DDMs, it is bounded by the maximum number of subdomains a collocation point can belong to, which is 4, and because of the least-squares problem, we obtain $4^2 = 16$. This results in a significantly lower condition number and GMRES will tend to converge more rapidly. 
	Meanwhile, the largest eigenvalue of the SAS preconditioned matrix is much smaller than that of the AS preconditioned matrix, and its distribution becomes more concentrated than that of RAS. This is because the RAS preconditioned matrix includes a Boolean matrix, which is introduced to reduce communication between different subdomains. However, a potential risk for the RaNN-based method is that this may increase linear dependencies between local basis functions, as zero values occur in their common region, which may lead to a higher condition number than SAS. Consequently, the SAS preconditioner, which avoids this issue, provides better conditioning and requires smaller iteration steps to converge.

	Besides, we try variants of the CG method to solve this problem. In Table \ref{table1b}, both CG and its variants work with the AS preconditioner, allowing them to utilize fewer number of iteration steps compared to GMRES. Additionally, the Conjugate Gradient Squared (CGS) and Biconjugate Gradient (BICG) methods are effective with the SAS and RAS preconditioners when CG fails. 
	
\begin{table}[!ht]	
	\centering  
	\renewcommand{\arraystretch}{1.1}
	\setlength\tabcolsep{4mm}
	\caption{{Relative $L^2$ error and number of iterations for different iterative solvers  in Example \ref{ex1}.}}	\label{table1b}	
	\begin{tabular}{ccccccccc}  
		\toprule 
		DoF & $N$ &  $\kappa(H^TH)$ & $M^{-1}$  & $\kappa(M^{-1}H^TH)$ & Solver  & iter  & $e_{L^2}$ \\ \midrule
		\multirow{4}{*}{256} & 	\multirow{4}{*}{1600} & \multirow{4}{*}{$10^{12}$ }  &  \multirow{4}{*}{$None$} & \multirow{4}{*}{$10^{12}$ }
		& CG & 2560 & 1.95e-2 \\
		&&&& & CGS & 2560 & 2.63e-2 \\
		&&&& & BICG & 2560 & 1.03e-2\\
		&&&& & GMRES & 2560  & 8.68e-2 \\
		\midrule
		\multirow{4}{*}{256} & 	\multirow{4}{*}{1600} & \multirow{4}{*}{$10^{12}$ }  &  \multirow{4}{*}{$M^{-1}_{AS}$} & \multirow{4}{*}{$10^{8}$}
		& CG & 8 & 5.03e-3 \\
		&&&& & CGS & 4 & 5.04e-3 \\
		&&&& & BICG & 8 & 5.08e-3\\
		&&&& & GMRES & 13 &  5.07e-3 \\
		\midrule
		\multirow{3}{*}{256} & 	\multirow{3}{*}{1600} & \multirow{3}{*}{$10^{12}$ }  &  \multirow{3}{*}{$M^{-1}_{SAS}$} & \multirow{3}{*}{$10^{8}$}
		& CGS & 6& 5.04e-3 \\
		&&&& & BICG & 11 & 5.09e-3\\
		&&&& & GMRES & 11&  5.08e-3 \\
		\midrule
		\multirow{3}{*}{256} & 	\multirow{3}{*}{1600} & \multirow{3}{*}{$10^{12}$ }  &  \multirow{3}{*}{$M^{-1}_{RAS}$} & \multirow{3}{*}{$10^{8}$}
		& CGS & 24 & 5.03e-3 \\
		&&&& & BICG & 32 & 5.05e-3\\
		&&&& & GMRES & 31 &  5.06e-3 \\
		
		\bottomrule
	\end{tabular}  
\end{table}

	To improve the accuracy, we can either increase the number of hidden nodes or use a larger overlap ratio. However, this also increases the condition number. For example, when setting $m = 32$ and $\delta = 2$, the matrix $H^TH$ has a large condition number around $10^{16}$, and the relative $L^2$ error can be around $10^{-5}$; when setting $m = 32$ and $\delta = 3$, the condition number increases further to approximately $10^{20}$, and the relative $L^2$ error  can be reduced to $10^{-7}$.  In this case, $H^TH$ would also be rank-deficient, and the preconditioned matrices would still exhibit a large condition number. By applying the PCA process, the condition number can then be significantly reduced, improving the convergence of the iterative methods.
	
			\begin{table}[!htbp]	
		\centering  
		\renewcommand{\arraystretch}{1.1}
		\setlength\tabcolsep{3mm}
		\caption{{Relative $L^2$ error and number of iterations for different preconditioners by GMRES in Example \ref{ex1}.}}	\label{table1c}	
		\begin{tabular}{ccccccccc}  
			\toprule 
			$\tau$ & DoF & $\kappa(H^TH)$ & $M^{-1}$  & $\kappa(M^{-1}H^TH)$ & $\sigma_{min}$ & $\sigma_{max}$ & iter  & $e_{L^2}$ \\ \midrule
			\multirow{3}{*}{$10^{-4}$} & 	\multirow{3}{*}{512} & \multirow{3}{*}{$10^{16}$}  & 
			$None$ & $ 10^{16}$ & $ 10^{-10}$& $10^{6}$ & 5120 & 3.72e-2\\
			&&& $M^{-1}_{AS}$ & $10^{12}$ & $10^{-6}$ & $10^{6}$ & 27 & 5.46e-5 \\
			&&& $M^{-1}_{SAS}$ & $ 10^{12}$  & $10^{-7}$ & $10^{5}$& 30 &  5.49e-5 \\
			\midrule
			\multirow{3}{*}{$\bf 10^{-3}$} & 	\multirow{3}{*}{\bf 436} & \multirow{3}{*}{$ \bf10^{13}$}  & 
			$None$ & $ 10^{13}$ & $10^{-8}$& $ 10^{5}$ & 4360& 3.75e-2\\
			&&& $\boldsymbol{M^{-1}_{AS}}$ & { $\bf 10^{10}$} & ${\bf 10^{-5}}$ & ${\bf 10^{5}}$& \textbf{16} &\textbf{1.28e-4} \\
			&&& $M^{-1}_{SAS}$ & $ 10^{10}$  & $10^{-6}$ & $10^{4}$& 18& 1.28e-4 \\
			\midrule
			\multirow{3}{*}{$10^{-2}$} & 	\multirow{3}{*}{335} & \multirow{3}{*}{$ 10^{10}$}  & 
			$None$ & $10^{10}$ & $10^{-5}$& $ 10^{5}$ & 3350 & 4.51e-2\\
			&&& $M^{-1}_{AS}$ & $ 10^{7}$ & $10^{-3}$ & $10^{4}$& 14& 7.14e-4 \\
			&&& $M^{-1}_{SAS}$ & $ 10^{7}$  & $10^{-4}$ & $10^{3}$&  13 &  7.11e-4 \\
			\midrule
			\multirow{3}{*}{$10^{-1}$} & 	\multirow{3}{*}{212} & \multirow{3}{*}{$ 10^{8}$}  & 
			$None$ & $10^{8}$ & $ 10^{-3}$& $ 10^{6}$ & 2120 & 5.01e-2\\
			&&& $M^{-1}_{AS}$ & $ 10^{4}$ & $10^{-2}$& $10^{3}$ & 12 & 7.13e-3 \\
			&&& $M^{-1}_{SAS}$ & $ 10^{4}$  &$10^{-3}$ &$10^{2}$& 11 & 7.10e-3 \\
			\bottomrule
		\end{tabular}  
	\end{table}
	
	In Table \ref{table1c}, let $m = 32$ and $\delta = 2$, we then examine how the threshold parameter $\tau$ for dropping small singular values based on the singular values influences the condition number and accuracy. It is varied from $10^{-4}$ to $10^{-1}$.  Within each subdomain $\Omega_j$, we apply the same threshold parameter $\tau$ to determine the number of effective neurons $p_j$. Then, the number of parameters that need to be solved is given by DoF $ = \sum\limits_{j=1}^{J} p_j$. We use the GMRES method without a preconditioner, and compare with the AS and SAS preconditioners. We observe that,  for both the AS and the SAS preconditioner, as the DoF decreases, the minimum singular value of  preconditioned matrices grows. This  leads to a lower condition number, which enables convergence within fewer iteration steps.
	
\begin{figure}[!ht]
	\includegraphics[width=6.5in]{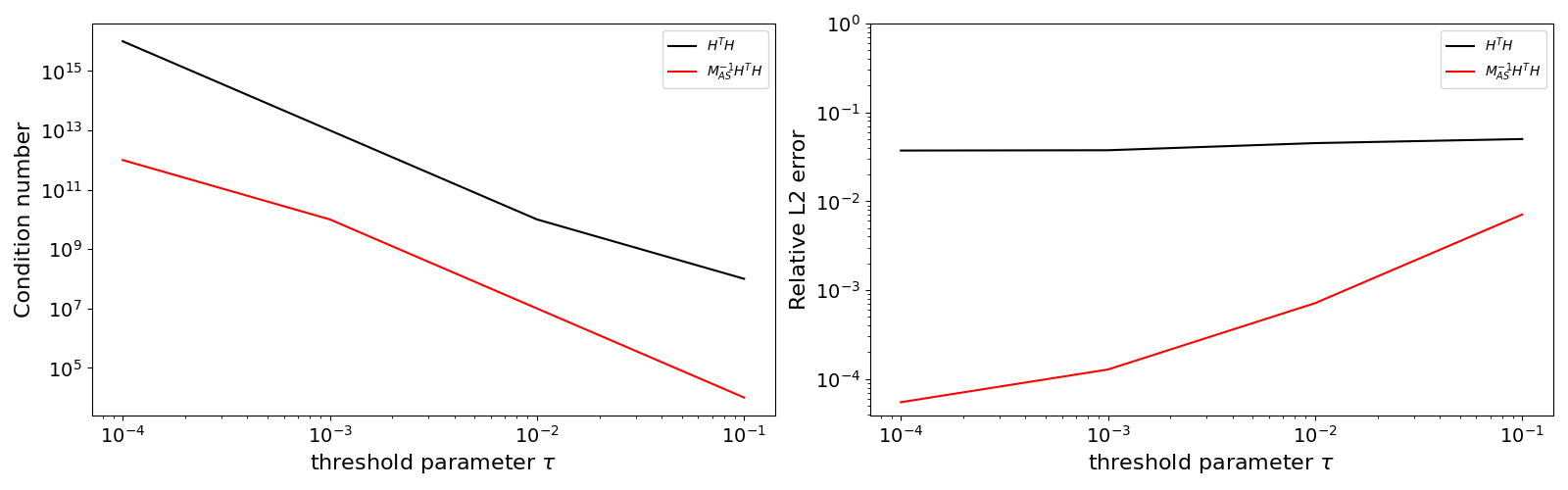}	
	\caption{The progressive of condition number and relative $L^2$ error via the  threshold parameter $\tau$   by different preconditioners in Example \ref{ex1}. }
	\label{figure3}
	\end{figure}

	Additionally, we observe that setting $\tau = 10^{-3}$ helps to lower the condition number without significant accuracy loss. Figure \ref{figure3} illustrates the progression of the condition number and relative $L^2$ error with respect to $\tau$, showing a substantial decrease in the condition number. The decrease is approximately two orders of magnitude, independent of $\tau$. Meanwhile, the relative $L^2$ error increases gradually, and does not exceed the result obtained by solving the problem without a preconditioner.

\begin{figure}[!ht] 		
		\centering
		\includegraphics[scale=0.14]{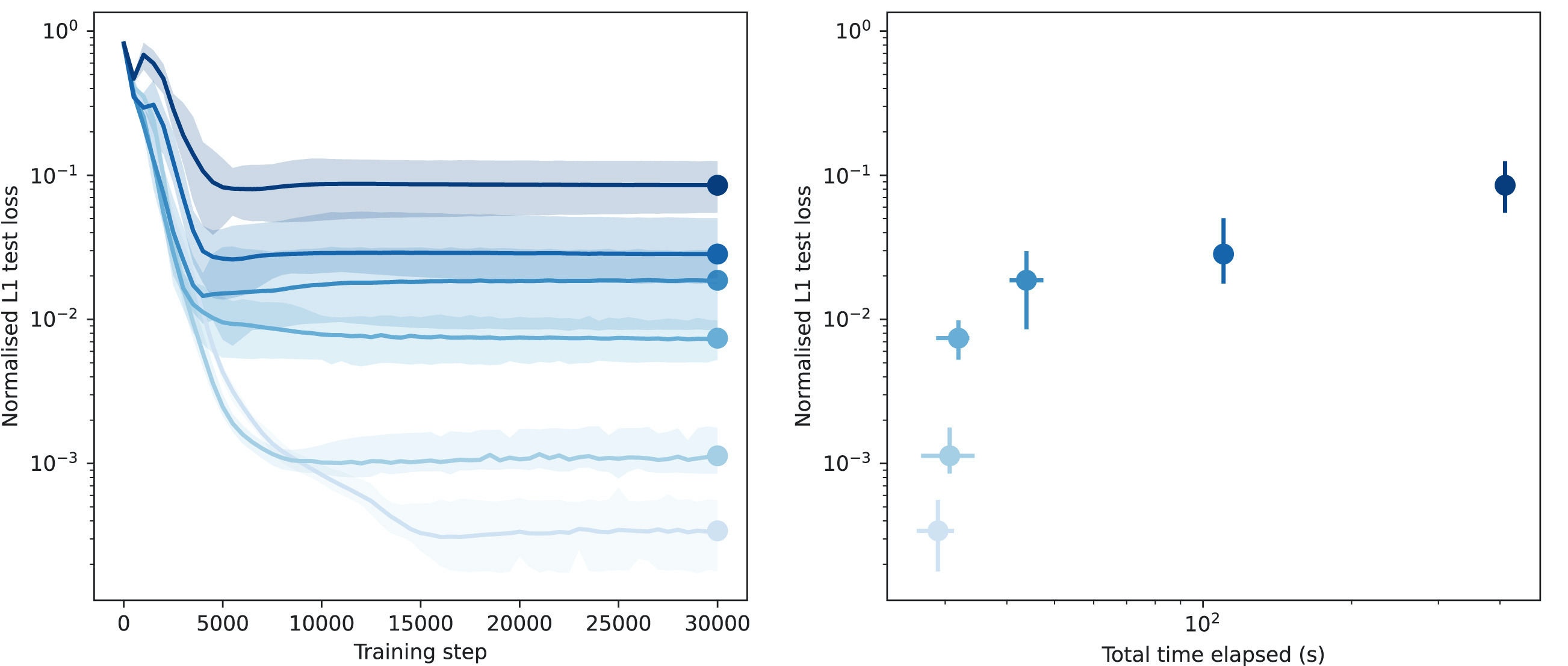}
		\caption{Weak scaling test using multilevel FBPINNs (taken from \cite{dolean2024multilevel}) in Example \ref{ex1}. The color-coded convergence curves and training times for each model are shown, and the color transitions from light to dark represent increasing problem complexity, ranging from  $n=1$ to $n=6$.}
		\label{Fbpinn}
	\end{figure}

	Finally, we consider using an increasing number of subdomains $J = 2^n \times 2^n$ to solve the problem \eqref{exlapla} with higher complexity $n$, and compare against the results obtained from multilevel FBPINNs \cite{dolean2024multilevel} with the weak scaling test shown in Figure \ref{Fbpinn}, which tested the same model problem. RaNNs with one hidden layer and 32 hidden nodes on each subdomain are employed. The number of collocation points is given to $N_{int} = (5 \times 2^{n}) \times (5 \times 2^{n})$, and the overlap ratio is fixed by  $\delta = 2$. The threshold parameter $\tau$ is fixed by $10^{-3}$.	The normalized $L^1$ test loss $e_{L^1}$ is computed approximately on $M = 350 \times 350$ uniformly-spaced test points by $\frac{1}{M} \sum_{i=1}^{M} \Vert \hat{u}(x_i) - u(x_i) \Vert_{l^1} / \gamma$, where $\gamma$ is defined as the standard deviation of the true solutions set $\{u(x_i)\}_{i=1}^{M}$. We use the normalized $L^1$ test loss here to compare with the results in \cite{dolean2024multilevel}.
			
	\begin{table}[!ht]	
		\centering  
		\renewcommand{\arraystretch}{1.1}
		\setlength\tabcolsep{3mm}
		\caption{{Weak scaling test using GMRES method for different preconditioners in Example \ref{ex1}.}}	\label{table2}	
		\begin{tabular}{cccccccccc}  
			\toprule 
			$n$ & $J$ & DoF & $N$ & $M^{-1}$  & $\kappa(M^{-1}H^TH)$ & iter  & $e_{L^1}$ & Time(s)\\ \midrule
			\multirow{3}{*}{2}  & \multirow{3}{*}{$2\times 2$} &  	 \multirow{3}{*}{109} &\multirow{3}{*}{$20 \times 20$}  & 
			$None$ & $10^{13}$ & 1090 & 6.81e-2 & 0.064\\
			&&&& $M^{-1}_{AS}$ & $1$ & 1& 1.98e-3  & 0.00040 \\
			&&&& $M^{-1}_{SAS}$ & $1$  & 1 & 1.97e-3 &   0.00033 \\
			\midrule
			\multirow{3}{*}{3}  & \multirow{3}{*}{$4 \times 4$} &   \multirow{3}{*}{436} & \multirow{3}{*}{$40 \times 40$}  & 
			$None$ & $10^{14}$  & 4330& 2.04e-1 & 0.81 \\
			&&&& $M^{-1}_{AS}$ & $ 10^{9}$ & 12 & 1.28e-2 & 0.00056 \\
			&&&& $M^{-1}_{SAS}$ & $ 10^{9}$  & 21 & 1.27e-2 &   0.00082 \\
			\midrule
			\multirow{3}{*}{4}  & \multirow{3}{*}{$8 \times 8$} &  \multirow{3}{*}{1706} &	\multirow{3}{*}{$80 \times 80$}  & 
			$None$ & $ 10^{14}$  & 17060 & 1.92e-1& 7.14 \\
			&&&& $M^{-1}_{AS}$ & $10^{10}$ & 38 &  2.46e-2  & 0.042\\
			&&&& $M^{-1}_{SAS}$ & $10^{10}$ & 52 & 2.48e-2  &   0.056 \\
			\midrule
			\multirow{3}{*}{5}  & \multirow{3}{*}{$16\times 16$} &  \multirow{3}{*}{6811} &	\multirow{3}{*}{$160 \times 160$}  & 
			$None$ & $10^{14}$ & 68110 & 3.24e-1 & 116 \\
			&&&& $M^{-1}_{AS}$ & $10^{10}$ & 51 & 3.24e-2 &  0.31\\
			&&&& $M^{-1}_{SAS}$ & $ 10^{10}$& 60 &    3.22e-2 & 0.35\\
			\midrule
			\multirow{3}{*}{\bf 6}  & \multirow{3}{*}{ $\bf32 \boldsymbol{\times} 32$} &  \multirow{3}{*}{\bf 27147} &	\multirow{3}{*}{$\bf320 \boldsymbol{\times} 320$}  & 
			$None$ & $10^{14}$ & 271470 & 4.36e-1 & 2950  \\
			&&&&  $\boldsymbol{M^{-1}_{AS}}$ & ${\bf 10^{11}}$ &  {\bf 109} & {\bf 6.71e-2} & {\bf 5.35} \\
			&&&& $M^{-1}_{SAS}$ & $10^{11}$& 138&   6.50e-2&5.67\\
			\bottomrule
		\end{tabular}  
	\end{table}
	
\begin{figure}[!htbp] 		
		\centering
		\includegraphics[scale=0.42]{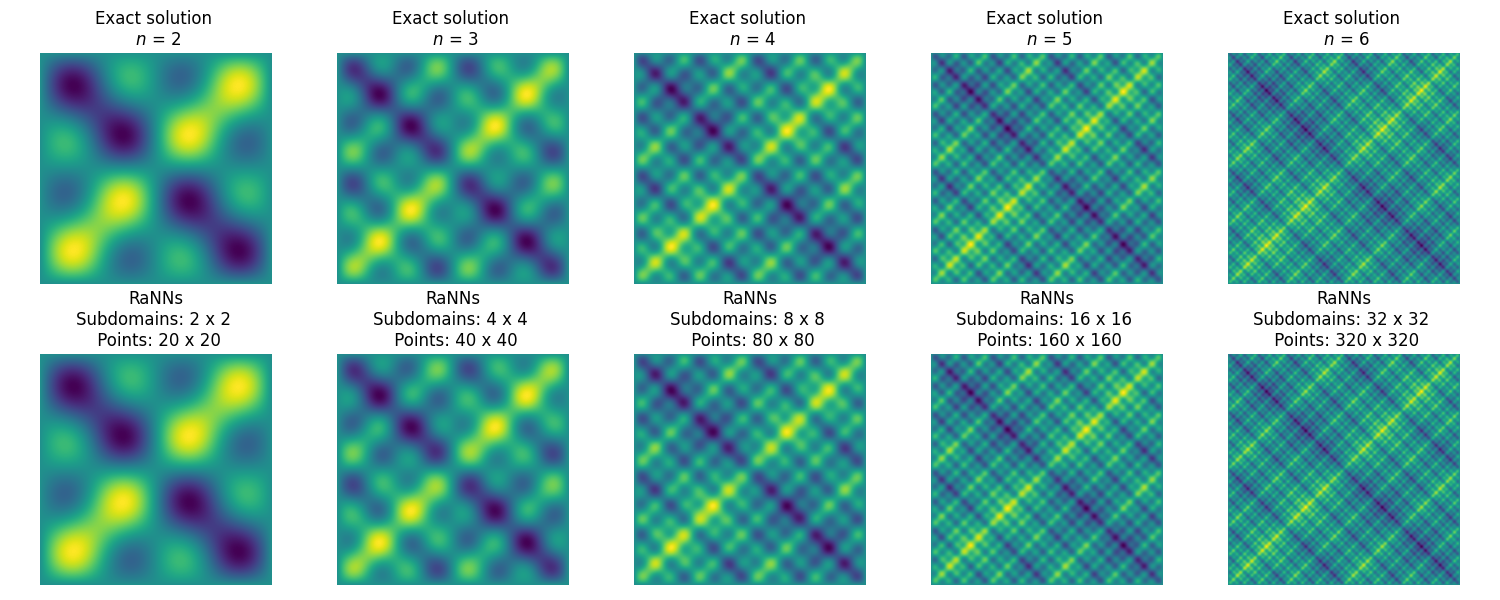}
		\caption{Weak scaling test using GMRES method with AS preconditioner in Example \ref{ex1}. The top figure shows the exact solution for different problem complexity $n$, while the bottom figure presents the numerical results corresponding to each $n$.}
		\label{Ranns}
	\end{figure}

	Table \ref{table2} shows the normalized $L^1$ test loss and the time for solving the preconditioned system using the GMRES method. We observe that both the AS and SAS preconditioners significantly reduce the condition number and result in fewer iterations, with the solution time not increasing drastically as the complexity $n$ increases. However, the normalized $L^1$ test loss does reduce slightly as the problem complexity increases. The exact and numerical solutions for $n=2,3,\ldots,6$ are shown in Figure \ref{Ranns}. Compared to the weak scaling test shown in Figure \ref{Fbpinn} for the multilevel FBPINNs method, our approach achieves better accuracy with a smaller domain decomposition, Specifically, only $32 \times 32$ subdomains are required to handle a problem with complexity $n = 6$. Meanwhile, the  time for solving the parameters of RaNNs using the GMRES method is much shorter than the training time for feedforward fully connected networks in multilevel FBPINNs.
	
	\vspace{-2mm}

	\begin{figure}[!ht] 		
		\centering
		\includegraphics[scale=0.4]{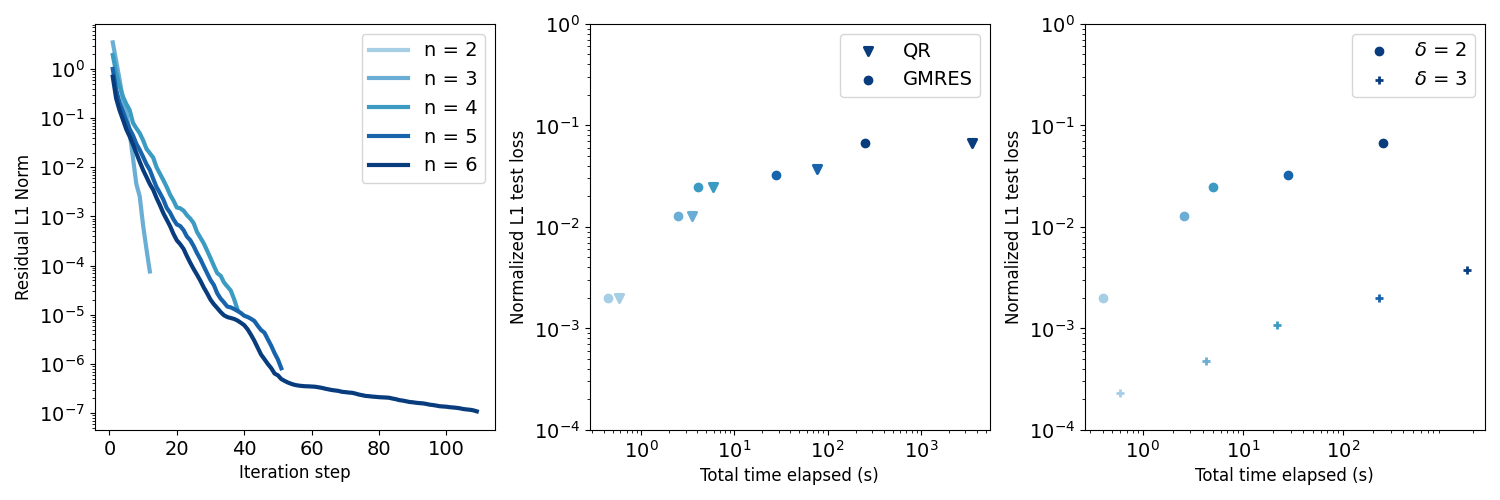}		
		\caption{Weak scaling test for the multi-scale Laplacian problem in Example \ref{ex1}. The left figure shows the progression of the residual $L^1$ norm with each iteration step when using the GMRES method with the AS preconditioner. The middle figure compares the total computation time and test loss for the QR decomposition method and GMRES method with the AS preconditioner. The right figure displays the total computation time for a different overlap ratio with $\delta = 2, 3$.}
		\label{ranns}
		
	\end{figure}

	Figure \ref{ranns} presents results from a weak scaling test using the AS preconditioner in Example \ref{ex1}. In the left plot, we observe the progression of the residual $\Vert H^THW - H^TF \Vert_{L^1}$ at each GMRES iteration for the problem \eqref{exlapla} with $n=2,3,\ldots,6$, illustrating the convergence behavior. The middle plot compares the  GMRES result against the results obtained using QR decomposition to directly solve the least-squares problem. Our approach maintains comparable accuracy to solving the least-squares problem using QR decomposition while substantially reducing computational time.  The right plot displays the results with a larger overlap ratio $\delta = 3$, where the normalized $L^1$ test loss is significantly reduced to approximately $10^{-3}$ without causing a major increase in computation time for each case. The total computational time includes the cost of constructing the linear system, building the preconditioner, and solving the preconditioned system using the GMRES method, with their respective percentages being approximately 60\%, 38\%, and 2\% for $n =6$. The primary part (linear system construction) will be further optimized in the future through parallelization.

	\begin{example}\label{exam2} Given $\Omega = (-1,1)$, we consider a one-dimensional advection-diffusion equation
	\begin{align*}
		\frac{\partial u}{\partial t}+	\frac{\partial u}{\partial x} &= \kappa \frac{\partial u^2}{\partial x^2} \qquad\qquad\quad \,{\rm in}\; \Omega\times I,  \\
		u(-1,t) &= u(1,t)  = 0 \quad \;\;\; \quad\; {\rm in}\;  I,\\
		u(x,0) &= -\sin(\pi x)\quad \; \qquad\,{\rm in}\; \Omega,
	\end{align*}
	with $I = (0,1)$. Here, $\kappa =0.1/\pi$ represents the diffusivity coefficient, which complicates the solution near the right boundary $x =1$ as the no-slip boundary condition is imposed \cite{kharazmi2021hp}.
\end{example}

We consider this equation as a 2-dimensional problem with a 1-dimensional spatial variable and a 1-dimensional temporal variable and solve it in the space-time domain instead of discretizing the temporal variable using the  finite difference method. We use the  Schwarz domain decomposition with $J = 4 \times 4 $ subdomains and compare the results obtained with the $hp$-VPINNs \cite{kharazmi2021hp}. The overlap ratio $\delta$ is fixed as 2. In \cite{kharazmi2021hp}, they employ individual networks in each subdomain with four layers and each with five neurons, resulting in 1296 training parameters, and a $4 \times 4 $ domain decomposition was adopted. RaNNs with one hidden layer and 81 hidden nodes on each subdomain are used to keep the same total number of parameters, and a uniform distribution $\mathcal{U}(-1,1)$ is applied to weights and bias parameters. The number of collocation points is set to $160 \times 160$, and the threshold parameter $\tau$ is fixed by $10^{-3}$. The constraining operator is chosen as 
\begin{equation} \label{cons2}
\mathcal{C}u(x) = L(x) u(x)+ G(x), \quad \text{where} \;	L(x)= \tanh(1+x) \tanh(1-x) \tanh(t), \; G(x) = -\sin(\pi x),
\end{equation}such that the boundary condition and initial condition are satisfied.

\begin{figure}[!ht]
	\begin{center}
		\subfigure[$hp$-VPINNs (taken from \cite{kharazmi2021hp})]{
			\centering
			\includegraphics[width=2.04in]{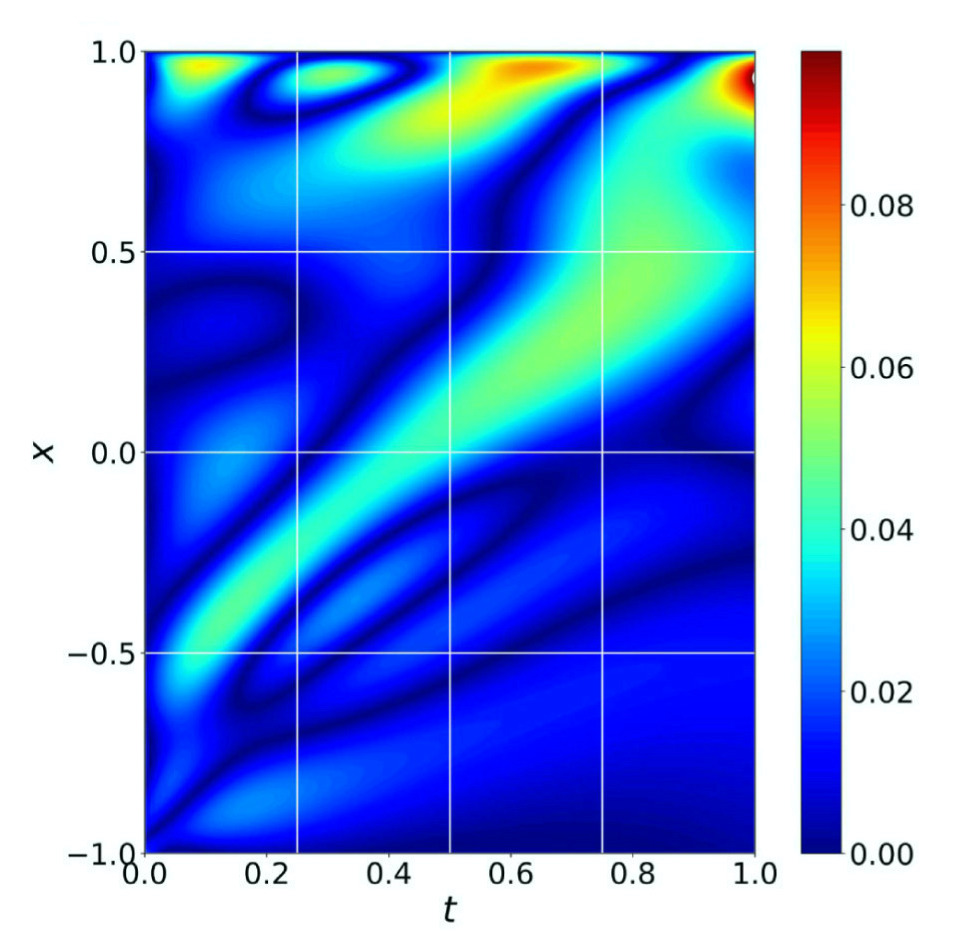}
		}
		\subfigure[Penalty method]{
			\centering
			\includegraphics[width=2.0in]{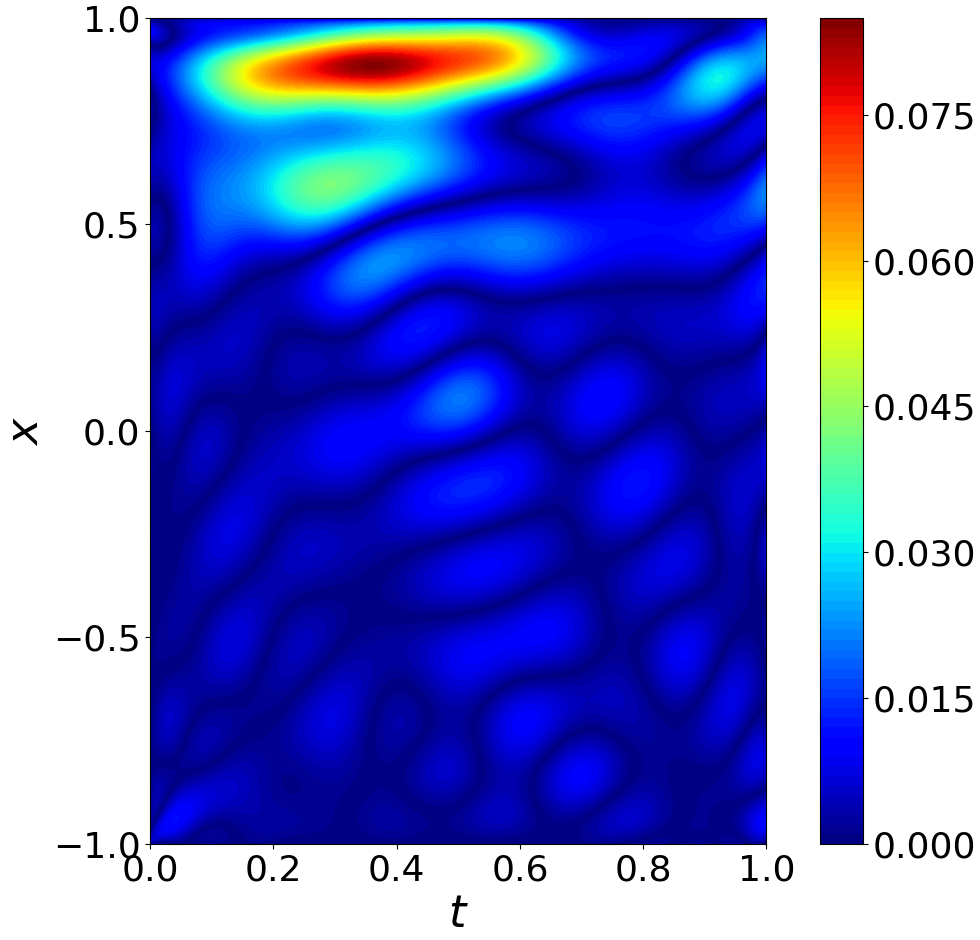}
		}
		\subfigure[Constraining operator]{
			\centering
			\includegraphics[width=2.0in]{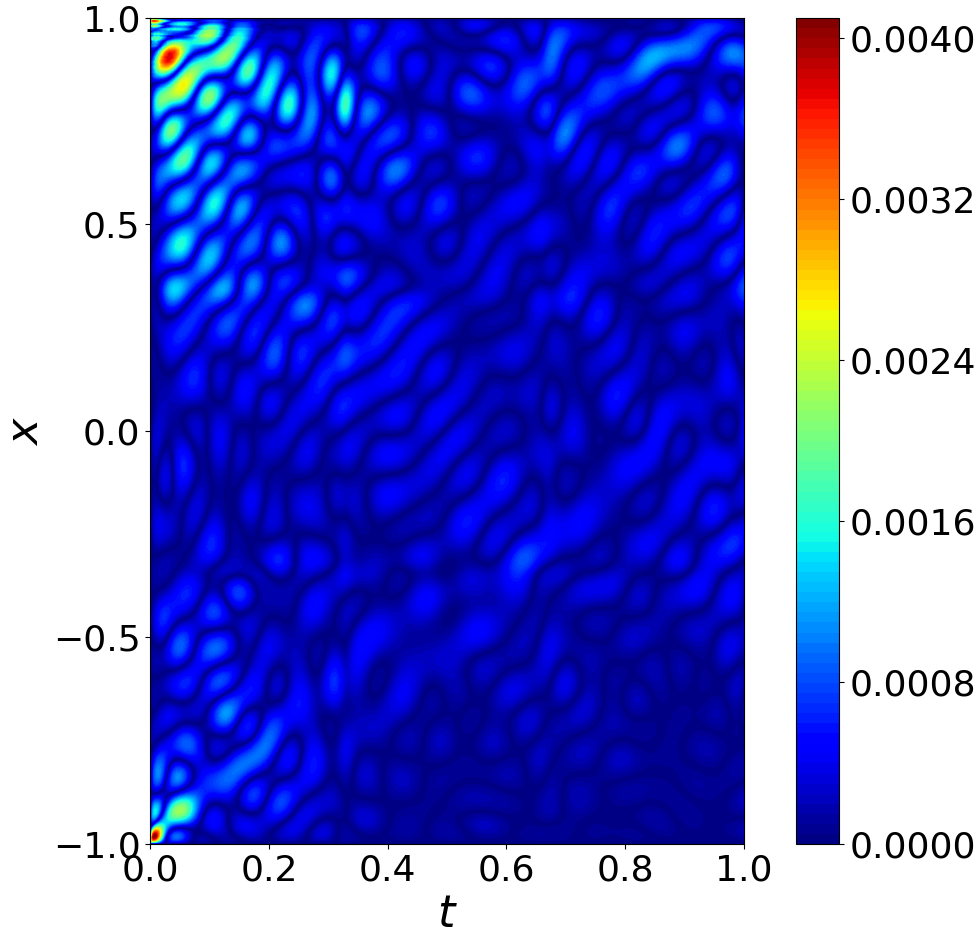}
		}
	\end{center}
	\vspace*{-15pt}
	\caption{Pointwise error of different methods with $ 4 \times 4 $ domain decomposition in Example \ref{exam2}. }
	\label{figure4}
\end{figure}

Figure \ref{figure4} illustrates the differences between the exact and numerical solutions obtained using various methods. The reference solution is calculated using a series summation with 800 terms of the solution in \cite{mojtabi2015one}. Figure \ref{figure4}(a) presents results computed using $hp$-VPINNs \cite{kharazmi2021hp}, while Figure \ref{figure4}(b) shows results using the penalty method to enforce boundary and initial conditions with a penalty parameter of 100. The computation takes 0.12 seconds and yields a relative $L^2$ error of $2.88\times10^{-2}$.  Figure \ref{figure4}(c) displays the results of applying the constraining operator computed via the CG method with an AS preconditioner, completing in 0.08 seconds after 138 iterations, and achieving a  relative $L^2$ error of $7.40\times10^{-4}$, around two orders of magnitude more accurate than $hp$-VPINNs and the penalty method. This demonstrates that our approach not only achieves higher accuracy and efficiency but also satisfies boundary and initial conditions directly, eliminating the need for penalty parameters.

\begin{example}\label{ex3} Given $\Omega = (0,1)^3$,  we consider a three-dimensional problem defined as follows:
	\begin{align*}
		  -\Delta \bu + \bu  &= \boldsymbol{f}  \quad\quad {\rm in}\; \Omega, \notag \\
		\bu(\bx) &= \boldsymbol{g} \;\;\;\;\;\;\;\, {\rm on}\; \partial \Omega. \notag 
	\end{align*}
	where the exact solution is given by:
	\begin{equation*}\bu = \begin{pmatrix}
			\cos(2\pi x_1)\sin(2\pi x_2)\sin(2\pi x_3) \\
			\sin(2\pi x_1)\cos(2\pi x_2)\sin(2\pi x_3)\\
			\sin(2\pi x_1)\sin(2\pi x_2)\cos(2\pi x_3)
		\end{pmatrix}
	\end{equation*}
The Dirichlet boundary condition $\boldsymbol{g}$ and the source term $\boldsymbol{f}$ are derived directly from the exact solution.
\end{example}

For this problem, we use RaNNs with one hidden layer and three output neurons to approximate $\bu = (u_1, u_2, u_3)$. The weights and bias parameters of the network are initialized using a uniform distribution $\mathcal{U}(-1,1)$. The $e_{L^2}$ error is computed approximately on $M = 100 \times 100 \times 100$ uniformly spaced test points. The overlap ratio $\delta$ is fixed as 2, and the threshold parameter $\tau$ is fixed as $10^{-3}$.

The constraining operator is chosen as $ \mathcal{C}\bu(\bx) = \boldsymbol{L} \bu(\bx)+ \boldsymbol{G}(\bx)$, where
\begin{align*} \label{cons3}
&\boldsymbol{L} (\bx)= \begin{pmatrix}
		L(\bx) \\
		L(\bx)  \\
		L(\bx)   \\	
	\end{pmatrix}, \; \boldsymbol{G}(x) = \begin{pmatrix}
	\sin(2\pi x_2)\sin(2\pi x_3) \\
	\sin(2\pi x_1)\sin(2\pi x_3)\\
	\sin(2\pi x_1)\sin(2\pi x_2)
\end{pmatrix}, \;\; \text{and}
	\\\;&{L} (\bx)  = \tanh(x_1)\tanh(1-x_1)\tanh(x_2)\tanh(1-x_2)\tanh(x_3)\tanh(1-x_3).
\end{align*}

In Table \ref{table3}, the domain is divided into $J = 2 \times 2 \times 2$ and $J = 4 \times 4 \times 4$ subdomains, respectively, while the number of hidden neurons $m$ is varied as 20, 40, and 80. We compare the results using the direct least-squares solver via QR decomposition, using the CG method with no preconditioner, and using the CG  method with the AS preconditioner. The number of iterations for the CG method is denoted by $iter$; we indicate 0 that represents the QR decomposition method, which does not require iterations. 

Table \ref{table3} illustrates that increasing the number of hidden neurons $m$ and subdomains improves the accuracy for both QR decomposition and CG method with an AS preconditioner. However, using the CG method  without a preconditioner may require a large number of iteration steps and often fails to converge to a satisfactory result. Meanwhile, the computational time required for the CG method  with the AS preconditioner is significantly lower than that for QR decomposition, particularly as the matrix size increases. This highlights the efficiency of the CG method with DD preconditioning in tackling 3-dimensional systems while maintaining accuracy.

\begin{table}[!htbp]	
	\centering  
	\renewcommand{\arraystretch}{1.1}
	\setlength\tabcolsep{3mm}
	\caption{{ Relative $L^2$ error and solving time for different methods in Example \ref{ex3}.}}	\label{table3}	
	\begin{tabular}{cccccccccc}  
		\toprule 
		$m$ & $J$ & DoF & $N$ & $M^{-1}$  & $\kappa(M^{-1}H^TH)$ & iter  & $e_{L^1}$ & Time(s)\\ \midrule
		\multirow{3}{*}{20}  & \multirow{3}{*}{$2\times 2 \times 2$} &  	 \multirow{3}{*}{480} &\multirow{3}{*}{$20 \times 20 \times 20$}  & 
		$None$ & $10^{7}$ & {0}& 1.35e-2 & 0.41 \\
		&&&&$None$ & $10^{7}$ & 1151 & 1.63e-2& 0.11\\
		&&&& $M^{-1}_{AS}$ & $10^{1}$ & 1 &1.35e-2&  0.001\\
		\midrule
		\multirow{3}{*}{40}  & \multirow{3}{*}{$2\times 2 \times 2$} &  	 \multirow{3}{*}{960} &\multirow{3}{*}{$20 \times 20 \times 20$}  & 
		$None$ & $10^{11}$ & {0}& 1.23e-3 & 1.04 \\
		&&&&$None$ & $10^{11}$ & 1175 & 2.98e-3& 0.43\\
		&&&& $M^{-1}_{AS}$ & $10^{1}$ & 1 &1.03e-3&  0.003\\
		\midrule
		\multirow{3}{*}{40}  & \multirow{3}{*}{$4\times 4 \times 4$} &  	 \multirow{3}{*}{7680} &\multirow{3}{*}{$40 \times 40 \times 40$}  & 
		$None$ & $10^{12}$ & {0}& 8.26e-4 & 181 \\
		&&&&$None$ & $10^{12}$ & 15236 & 2.38e-3& 385\\
		&&&& $M^{-1}_{AS}$ & $10^{7}$ & 10& 8.26e-4&  0.41 \\
		\midrule
		\multirow{3}{*}{80}  & \multirow{3}{*}{$4\times 4 \times 4$} &  	 \multirow{3}{*}{15360} &\multirow{3}{*}{$40 \times 40 \times 40$}  & 
		$None$ & $10^{15}$ & {0}& 2.37e-5 & 958 \\
		&&&&$None$ & $10^{15}$ &15984& 3.02e-4& 1069\\
		&&&& $M^{-1}_{AS}$ & $10^{10}$ & 11&  2.34e-5&  2.21 \\
		\bottomrule
	\end{tabular}  
\end{table}

\begin{figure}[!ht]
	\includegraphics[width=6.5in]{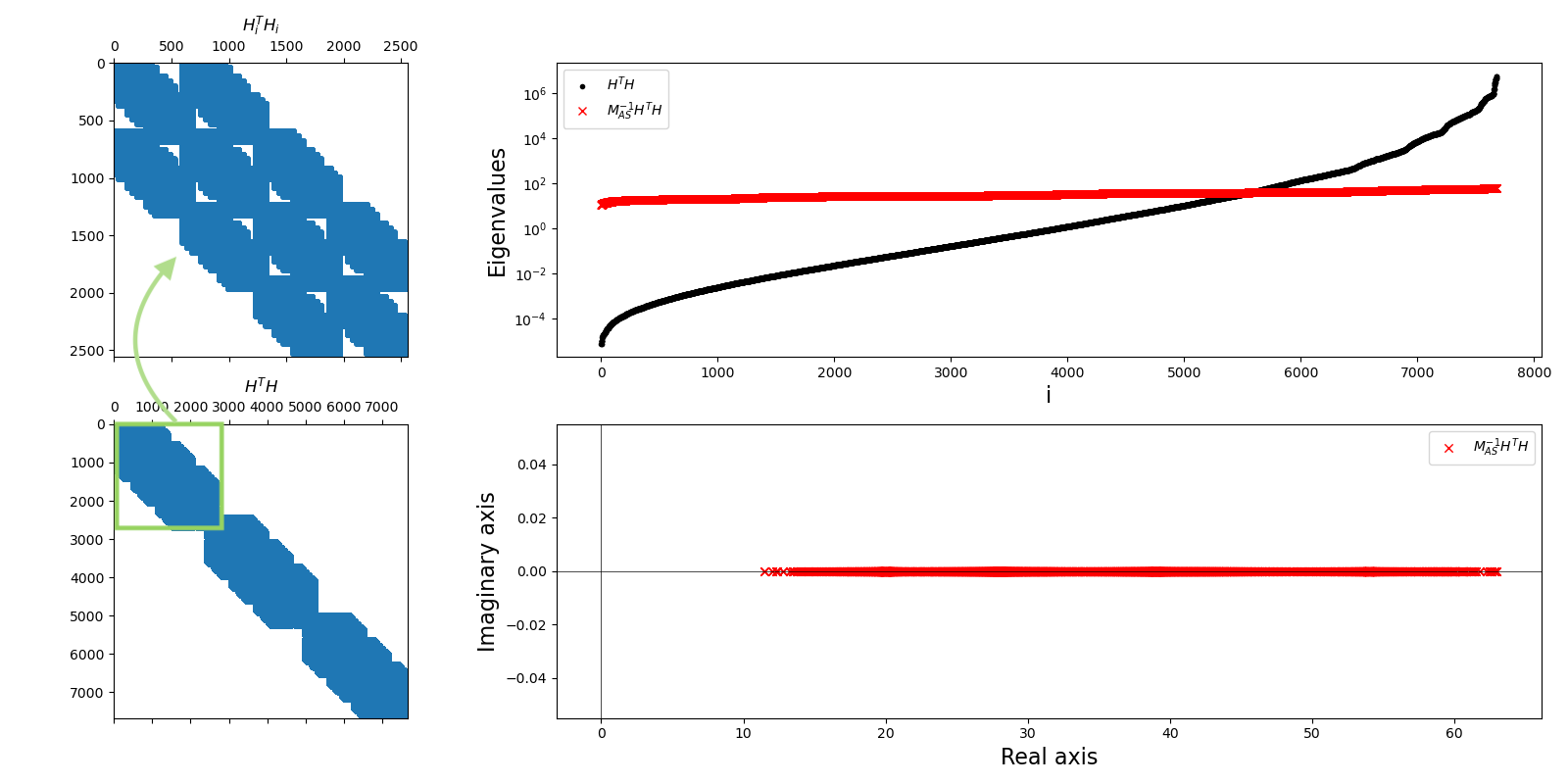}	
	\caption{Matrix information and the distribution of eigenvalues of the preconditioned matrix with the AS preconditioner for a $4 \times 4 \times 4$ domain decomposition in Example \ref{ex3}. The left figure presents the matrix information for $H^TH$ and its diagonal part. The top right figure compares the distribution of eigenvalues for the AS preconditioner against the distribution without preconditioning. The bottom right figure illustrates the eigenvalue distributions for the AS preconditioner in the complex plane.}
	\label{figure5}
\end{figure}

As shown on the left of Figure \ref{figure5}, we display the matrix information for a $4 \times 4 \times 4$ domain decomposition with $m = 40$ in Example \ref{ex3}. To construct the preconditioner for this equation, we simply need to diagonally stack the preconditioners constructed for $H_i^TH_i$, which corresponds to a scalar problem of the form
\begin{align*}
	-\Delta u_i + u_i  &= {f_i}  \quad\quad {\rm in}\; \Omega, \; i = 1,2,3, \notag \\
	u_i(\bx) &= {g_i} \;\;\quad\,\; {\rm on}\; \partial \Omega, \; i = 1,2,3. \notag 
\end{align*}
 The overall preconditioner remains a block diagonal matrix. In Figure \ref{figure5}, the right figure displays the distribution of all eigenvalues for $H^T H$ and $M^{-1}_{AS} H^T H$, while the bottom shows the eigenvalue distribution in the complex plane for the AS preconditioner. We observe that the eigenvalue distribution of the AS preconditioned matrix becomes more concentrated, effectively shifting the minimum eigenvalue from $10^{-6}$ to $11$, with the maximum eigenvalue being smaller than 64. In classical DDMs, it is bounded by  the maximum number of subdomains a collocation point can belong to, which is 8. Due to the least-squares problem, this number has to be squared, resulting in $8^2 = 64$. This results in a more concentrated distribution and a significantly  lower condition number.

\section{Conclusion}

In this research, we present a new approach for solving PDEs using randomized neural networks  across overlapping subdomains interconnected by local window functions. By integrating a PCA scheme into the neural network structure on the local subdomains, we significantly reduce the number of parameters by eliminating less essential information of neural networks, resulting in a system with a lower condition number. Furthermore, constructing overlapping Schwarz preconditioners for the resulting system enables the efficient solution of the least-squares problem using preconditioned iterative solvers by further reducing the condition number. This means that overlapping Schwarz domain decomposition is combined with RaNNs in two main ways: first, to formulate the least-squares problem, and second, to efficiently construct preconditioners for the resulting linear systems. This is the first paper to demonstrate the application of overlapping Schwarz preconditioners for RaNNs with domain decomposition. The main advantage of this method is the significant reduction in training time compared to traditional neural network training techniques without loss of accuracy.

The integration of RaNNs with DDMs demonstrates considerable potential for addressing multi-scale problems and time-dependent PDEs, offering a robust and efficient solution. However, to solve more complex problems, further investigation into RaNNs combined with multi-level DDM-based architectures is necessary. Specifically, the development of two-level or multi-level Schwarz preconditioners should be pursued. To manage the increased complexity, parallel methods must be developed to ensure computational efficiency. Moreover, careful construction of the corresponding preconditioners is essential, particularly as the system may include dense blocks. Future research should also focus on the numerical analysis of this approach to further refine its theoretical foundation and practical applicability.

\bigskip
\noindent{\bf Data Availability Statement.} 
All data, models, or code that support the findings of this study are available from the corresponding author upon reasonable request.

\bibliographystyle{abbrv}
\bibliography{reference}

\begin{thebibliography}{10}

\bibitem{abdi2010principal}
H.~Abdi and L.~J. Williams.
\newblock Principal component analysis.
\newblock {\em Wiley interdisciplinary reviews: computational statistics},
  2(4):433--459, 2010.

\bibitem{abid2002new}
S.~Abid, F.~Fnaiech, and M.~Najim.
\newblock A new neural network pruning method based on the singular value
  decomposition and the weight initialisation.
\newblock In {\em 2002 11th European Signal Processing Conference}, pages 1--4.
  IEEE, 2002.

\bibitem{anderson2024elm}
S.~Anderson, V.~Dolean, B.~Moseley, and J.~Pestana.
\newblock {ELM-FBPINN}: efficient finite-basis physics-informed neural
  networks.
\newblock {\em arXiv preprint arXiv:2409.01949}, 2024.

\bibitem{cai1999restricted}
X.-C. Cai and M.~Sarkis.
\newblock A restricted additive {Schwarz} preconditioner for general sparse
  linear systems.
\newblock {\em Siam journal on scientific computing}, 21(2):792--797, 1999.

\bibitem{chan1994domain}
T.~F. Chan and T.~P. Mathew.
\newblock Domain decomposition algorithms.
\newblock {\em Acta numerica}, 3:61--143, 1994.

\bibitem{chen2022bridging}
J.~Chen, X.~Chi, W.~E, and Z.~Yang.
\newblock Bridging traditional and machine learning-based algorithms for
  solving pdes: the random feature method.
\newblock {\em J Mach Learn}, 1:268--98, 2022.

\bibitem{chen2024optimization}
J.~Chen, W.~E, and Y.~Sun.
\newblock Optimization of random feature method in the high-precision regime.
\newblock {\em Communications on Applied Mathematics and Computation},
  6(2):1490--1517, 2024.

\bibitem{Dang2024local}
H.~Dang and F.~Wang.
\newblock Local randomized neural networks with hybridized discontinuous
  {Petrov--Galerkin} methods for {Stokes--Darcy} flows.
\newblock {\em Physics of Fluids}, 36(8), 2024.

\bibitem{dang2024adaptive}
H.~Dang, F.~Wang, and S.~Jiang.
\newblock Adaptive growing randomized neural networks for solving partial
  differential equations.
\newblock {\em arXiv preprint arXiv:2408.17225}, 2024.

\bibitem{dissanayake1994neural}
M.~G. Dissanayake and N.~Phan-Thien.
\newblock Neural-network-based approximations for solving partial differential
  equations.
\newblock {\em Communications in Numerical Methods in Engineering},
  10(3):195--201, 1994.

\bibitem{dolean2022finite}
V.~Dolean, A.~Heinlein, S.~Mishra, and B.~Moseley.
\newblock Finite basis physics-informed neural networks as a {Schwarz} domain
  decomposition method.
\newblock In {\em International Conference on Domain Decomposition Methods},
  pages 165--172. Springer, 2022.

\bibitem{dolean2024multilevel}
V.~Dolean, A.~Heinlein, S.~Mishra, and B.~Moseley.
\newblock Multilevel domain decomposition-based architectures for
  physics-informed neural networks.
\newblock {\em Computer Methods in Applied Mechanics and Engineering},
  429:117116, 2024.

\bibitem{dolean2015introduction}
V.~Dolean, P.~Jolivet, and F.~Nataf.
\newblock {\em An introduction to domain decomposition methods: algorithms,
  theory, and parallel implementation}.
\newblock SIAM, 2015.

\bibitem{dong2021local}
S.~Dong and Z.~Li.
\newblock Local extreme learning machines and domain decomposition for solving
  linear and nonlinear partial differential equations.
\newblock {\em Computer Methods in Applied Mechanics and Engineering},
  387:114129, 2021.

\bibitem{gallicchio2020deep}
C.~Gallicchio and S.~Scardapane.
\newblock Deep randomized neural networks.
\newblock In {\em Recent Trends in Learning From Data: Tutorials from the INNS
  Big Data and Deep Learning Conference (INNSBDDL2019)}, pages 43--68.
  Springer, 2020.

\bibitem{hayashi1993fast}
M.~Hayashi.
\newblock A fast algorithm for the hidden units in a multilayer perceptron.
\newblock In {\em Proceedings of 1993 International Conference on Neural
  Networks (IJCNN-93-Nagoya, Japan)}, volume~1, pages 339--342. IEEE, 1993.

\bibitem{heinlein2024multifidelity}
A.~Heinlein, A.~A. Howard, D.~Beecroft, and P.~Stinis.
\newblock Multifidelity domain decomposition-based physics-informed neural
  networks for time-dependent problems.
\newblock {\em arXiv preprint arXiv:2401.07888}, 2024.

\bibitem{heinlein2021combining}
A.~Heinlein, A.~Klawonn, M.~Lanser, and J.~Weber.
\newblock Combining machine learning and domain decomposition methods for the
  solution of partial differential equations—a review.
\newblock {\em GAMM-Mitteilungen}, 44(1):e202100001, 2021.

\bibitem{howard2024finite}
A.~A. Howard, B.~Jacob, S.~H. Murphy, A.~Heinlein, and P.~Stinis.
\newblock Finite basis kolmogorov-arnold networks: domain decomposition for
  data-driven and physics-informed problems.
\newblock {\em arXiv preprint arXiv:2406.19662}, 2024.

\bibitem{huang2006extreme}
G.-B. Huang, Q.-Y. Zhu, and C.-K. Siew.
\newblock Extreme learning machine: theory and applications.
\newblock {\em Neurocomputing}, 70(1-3):489--501, 2006.

\bibitem{igelnik1995stochastic}
B.~Igelnik and Y.-H. Pao.
\newblock Stochastic choice of basis functions in adaptive function
  approximation and the functional-link net.
\newblock {\em IEEE transactions on Neural Networks}, 6(6):1320--1329, 1995.

\bibitem{jagtap2020extended}
A.~D. Jagtap and G.~E. Karniadakis.
\newblock Extended physics-informed neural networks ({XPINNs}): A generalized
  space-time domain decomposition based deep learning framework for nonlinear
  partial differential equations.
\newblock {\em Communications in Computational Physics}, 28(5), 2020.

\bibitem{jagtap2020conservative}
A.~D. Jagtap, E.~Kharazmi, and G.~E. Karniadakis.
\newblock Conservative physics-informed neural networks on discrete domains for
  conservation laws: Applications to forward and inverse problems.
\newblock {\em Computer Methods in Applied Mechanics and Engineering},
  365:113028, 2020.

\bibitem{kharazmi2021hp}
E.~Kharazmi, Z.~Zhang, and G.~E. Karniadakis.
\newblock {hp-VPINNs}: Variational physics-informed neural networks with domain
  decomposition.
\newblock {\em Computer Methods in Applied Mechanics and Engineering},
  374:113547, 2021.

\bibitem{klawonn2024machine}
A.~Klawonn, M.~Lanser, and J.~Weber.
\newblock Machine learning and domain decomposition methods-a survey.
\newblock {\em Computational Science and Engineering}, 1(1):2, 2024.

\bibitem{lagaris1998artificial}
I.~E. Lagaris, A.~Likas, and D.~I. Fotiadis.
\newblock Artificial neural networks for solving ordinary and partial
  differential equations.
\newblock {\em IEEE transactions on neural networks}, 9(5):987--1000, 1998.

\bibitem{leake2020deep}
C.~Leake and D.~Mortari.
\newblock Deep theory of functional connections: A new method for estimating
  the solutions of partial differential equations.
\newblock {\em Machine Learning and Knowledge Extraction}, 2(1):37--55, 2020.

\bibitem{li2020deep}
W.~Li, X.~Xiang, and Y.~Xu.
\newblock Deep domain decomposition method: Elliptic problems.
\newblock In {\em Mathematical and Scientific Machine Learning}, pages
  269--286. PMLR, 2020.

\bibitem{liu2014extreme}
X.~Liu, S.~Lin, J.~Fang, and Z.~Xu.
\newblock Is extreme learning machine feasible? a theoretical assessment (part
  i).
\newblock {\em IEEE Transactions on Neural Networks and Learning Systems},
  26(1):7--20, 2014.

\bibitem{luenberger1973introduction}
D.~G. Luenberger.
\newblock Introduction to linear and nonlinear programming.
\newblock {\em Addison-Wesley Reading, MA}, 28, 1973.

\bibitem{luo2021phase}
T.~Luo, Z.-Q.~J. Xu, Z.~Ma, and Y.~Zhang.
\newblock Phase diagram for two-layer {ReLU} neural networks at infinite-width
  limit.
\newblock {\em Journal of Machine Learning Research}, 22(71):1--47, 2021.

\bibitem{lyu2020enforcing}
L.~Lyu, K.~Wu, R.~Du, and J.~Chen.
\newblock Enforcing exact boundary and initial conditions in the deep mixed
  residual method.
\newblock {\em arXiv preprint arXiv:2008.01491}, 2020.

\bibitem{mcfall2009artificial}
K.~S. McFall and J.~R. Mahan.
\newblock Artificial neural network method for solution of boundary value
  problems with exact satisfaction of arbitrary boundary conditions.
\newblock {\em IEEE transactions on neural networks}, 20(8):1221--1233, 2009.

\bibitem{mojtabi2015one}
A.~Mojtabi and M.~O. Deville.
\newblock One-dimensional linear advection--diffusion equation: Analytical and
  finite element solutions.
\newblock {\em Computers \& Fluids}, 107:189--195, 2015.

\bibitem{moseley2023finite}
B.~Moseley, A.~Markham, and T.~Nissen-Meyer.
\newblock Finite basis physics-informed neural networks (fbpinns): a scalable
  domain decomposition approach for solving differential equations.
\newblock {\em Advances in Computational Mathematics}, 49(4):62, 2023.

\bibitem{pao1994learning}
Y.-H. Pao, G.-H. Park, and D.~J. Sobajic.
\newblock Learning and generalization characteristics of the random vector
  functional-link net.
\newblock {\em Neurocomputing}, 6(2):163--180, 1994.

\bibitem{pao1992functional}
Y.-H. Pao and Y.~Takefuji.
\newblock Functional-link net computing: theory, system architecture, and
  functionalities.
\newblock {\em Computer}, 25(5):76--79, 1992.

\bibitem{psichogios1994svd}
D.~C. Psichogios and L.~H. Ungar.
\newblock Svd-net: An algorithm that automatically selects network structure.
\newblock {\em IEEE Transactions on Neural Networks}, 5(3):513--515, 1994.

\bibitem{raissi2019physics}
M.~Raissi, P.~Perdikaris, and G.~E. Karniadakis.
\newblock Physics-informed neural networks: A deep learning framework for
  solving forward and inverse problems involving nonlinear partial differential
  equations.
\newblock {\em Journal of Computational physics}, 378:686--707, 2019.

\bibitem{schwarz1870ueber}
H.~A. Schwarz.
\newblock {\em Ueber einen Grenz{\"u}bergang durch alternirendes Verfahren}.
\newblock Z{\"u}rcher u. Furrer, 1870.

\bibitem{shang2023Randomized}
Y.~Shang, F.~Wang, and J.~Sun.
\newblock Randomized neural network with {Petrov-Galerkin} methods for solving
  linear and nonlinear partial differential equations.
\newblock {\em Communications in Nonlinear Science and Numerical Simulation},
  127(Dec.):1.1--1.20, 2023.

\bibitem{smithdomain}
B.~Smith, P.~Bj{\o}rstad, and W.~Gropp.
\newblock Domain decomposition, parallel multilevel methods for elliptic
  partial differential equations,(1996).

\bibitem{sun2024local}
J.~Sun, S.~Dong, and F.~Wang.
\newblock Local randomized neural networks with discontinuous {Galerkin}
  methods for partial differential equations.
\newblock {\em Journal of Computational and Applied Mathematics}, 445:115830,
  2024.

\bibitem{sun2024local_W}
J.~Sun and F.~Wang.
\newblock Local randomized neural networks with discontinuous {Galerkin}
  methods for diffusive-viscous wave equation.
\newblock {\em Computers \& Mathematics with Applications}, 154:128--137, 2024.

\bibitem{tamura1993determination}
S.~Tamura, M.~Tateishi, M.~Matumoto, and S.~Akita.
\newblock Determination of the number of redundant hidden units in a
  three-layered feedforward neural network.
\newblock In {\em Proceedings of 1993 International Conference on Neural
  Networks (IJCNN-93-Nagoya, Japan)}, volume~1, pages 335--338. IEEE, 1993.

\bibitem{teoh2006estimating}
E.~J. Teoh, K.~C. Tan, and C.~Xiang.
\newblock Estimating the number of hidden neurons in a feedforward network
  using the singular value decomposition.
\newblock {\em IEEE Transactions on Neural Networks}, 17(6):1623--1629, 2006.

\bibitem{toselli2006domain}
A.~Toselli and O.~Widlund.
\newblock {\em Domain decomposition methods-algorithms and theory}, volume~34.
\newblock Springer Science \& Business Media, 2006.

\bibitem{weigend1991effective}
A.~S. Weigend and D.~E. Rumelhart.
\newblock The effective dimension of the space of hidden units.
\newblock In {\em [Proceedings] 1991 IEEE International Joint Conference on
  Neural Networks}, pages 2069--2074. IEEE, 1991.

\bibitem{widlund1988iterative}
O.~B. Widlund.
\newblock Iterative substructuring methods: Algorithms and theory for elliptic
  problems in the plane.
\newblock In {\em First International Symposium on Domain Decomposition Methods
  for Partial Differential Equations, Philadelphia, PA}, pages 113--128, 1988.

\bibitem{xue2013restructuring}
J.~Xue, J.~Li, and Y.~Gong.
\newblock Restructuring of deep neural network acoustic models with singular
  value decomposition.
\newblock In {\em Interspeech}, pages 2365--2369, 2013.

\bibitem{zhang2024transferable}
Z.~Zhang, F.~Bao, L.~Ju, and G.~Zhang.
\newblock Transferable neural networks for partial differential equations.
\newblock {\em Journal of Scientific Computing}, 99(1):2, 2024.

\bibitem{zhou2022towards}
H.~Zhou, Z.~Qixuan, T.~Luo, Y.~Zhang, and Z.-Q. Xu.
\newblock Towards understanding the condensation of neural networks at initial
  training.
\newblock {\em Advances in Neural Information Processing Systems},
  35:2184--2196, 2022.

\end{thebibliography}

\end{document}